\documentclass[3p]{elsarticle}

\usepackage{etex} 
\reserveinserts{30}
\usepackage{amsmath, amssymb}
\usepackage{graphicx}
\usepackage{tikz}
\usepackage{xcolor,calc}
\usetikzlibrary{matrix,arrows,decorations.pathmorphing}
\usepackage{subfigure}
\usepackage{pgfplots}
\usepackage{amsmath,amscd}%
\usepackage{amsfonts}%
\usepackage{amssymb}%
\usepackage{graphicx}
\usepackage{txfonts}
\usepackage{color}
\usepackage[numbers]{natbib}

\newtheorem{example}{Example}

\newtheorem{remark}{Remark}

\numberwithin{equation}{section}

\newcommand{\kdifform}[2]{#1^{(#2)}} 
\newcommand{\kdifformh}[2]{#1^{(#2)}_{h}} 
\newcommand{\kchain}[2]{\mathbf{#1}_{(#2)}} 
\newcommand{\kcochain}[2]{\mathbf{#1}^{(#2)}} 
\newcommand{\ccomplex}[1]{{#1}} 
\newcommand{\reconstruction}{\mathcal{I}} 
\newcommand{\reduction}{\mathcal{R}} 
\newcommand{\kchainspacedomain}[2]{C_{#1}(#2)} 
\newcommand{\incidenceboundary}[2]{\mathsf{E}_{(#1,#2)}} 
\newcommand{\innerspace}[3]{\left(#1,#2\right)_{#3}} 
\newcommand{\matrixoperator}[1]{#1}
\newcommand{\inner}[2]{\left( #1, #2\right)} 

\renewcommand{\a}{\alpha}
\renewcommand{\b}{\beta}

\newcommand{\projection}{\pi_{h}} 
\newcommand{\dualprojection}{\tilde{\pi}_{h}} 
\newcommand{\pullback}{\Phi^{\star}} 

\newcommand{\manifold}[1]{\mathcal{#1}} 
\newcommand{\tangentspace}[1]{T_{p}\manifold{#1}} 
 
\newcommand{\ederiv}{\mathrm{d}} 
\newcommand{\coderiv}{\mathrm{d}^{*}} 

\newcommand{\duality}[2]{\langle #1, #2\rangle} 

\renewcommand{\eqref}[1]{(\ref{#1})} 

\newcommand{\realNumber}{\mathbb{R}}
\newcommand{\differentialFormSpace}[2]{\Lambda^{#1}\left(#2\right)}
\newcommand{\differentialFormSpaceDual}[2]{\tilde{\Lambda}^{#1}\left(#2 \right)}
\newcommand{\differentialFormSpaceDiscrete}[2]{\Lambda^{#1}_{h}\left(#2\right)}
\newcommand{\reductionDual}{\mathcal{\tilde{R}}}
\newcommand{\reconstructionDual}{\mathcal{\tilde{I}}}
\newcommand{\cochainSpace}[2]{ C^{#1}\left( #2 \right)}
\newcommand{\cochainSpaceDual}[2]{ \tilde{C}^{#1}\left( #2 \right)}

\journal{Journal of Computational Physics}

\begin{document}

\begin{frontmatter}


\author[tudelftw]{Artur Palha}
\ead{A.Palha@TUDelft.nl}

\author[tudelft]{Pedro Pinto Rebelo}
\ead{P.J.PintoRebelo@TUDelft.nl}

\author[tudelft]{Ren\'{e} Hiemstra}
\ead{R.R.Hiemstra@TUDelft.nl}

\author[shell]{Jasper Kreeft}
\ead{Jasper.Kreeft@Shell.com}

\author[tudelft]{Marc Gerritsma\corref{cor1}}
\ead{M.I.Gerritsma@TUDelft.nl}

\cortext[cor1]{Corresponding author}

\address[tudelftw]{Delft University of Technology, Faculty of Aerospace Engineering, Wind Energy Group P.O. Box 5058, 2600 GB Delft, The Netherlands}
\address[tudelft]{Delft University of Technology, Faculty of Aerospace Engineering, Aerodynamics Group P.O. Box 5058, 2600 GB Delft, The Netherlands}
\address[shell]{Shell Global Solutions, The Netherlands}


\title{Physics-compatible discretization techniques on single and dual grids, with application to the Poisson equation of volume forms}

\begin{abstract}
This paper introduces the basic concepts for physics-compatible discretization techniques. The paper gives a clear distinction between vectors and forms. Based on the difference between forms and pseudo-forms and the $\star$-operator which switches between the two, a dual grid description and a single grid description are presented. The dual grid method resembles a staggered finite volume method, whereas the single grid approach shows a strong resemblance with a finite element method. Both approaches are compared for the Poisson equation for volume forms.
\end{abstract}

\begin{keyword}
Mimetic discretization, differential forms, single grid, dual grid, geometric flexibility.
\end{keyword}

\end{frontmatter}

\section{INTRODUCTION}
\label{Section::Introduction}

Mimetic methods aim to preserve essential physical/mathematical structures in a discrete setting. Many of such structures are {\em topological}, i.e. independent of metric, and involve {\em integral relations}. Since integration will play an important role and integration of differential forms is a metric-free operation, we will work with differential forms. Formally, differentials forms are linear functionals on multi-vectors, but Flanders, \citep[p.1]{flanders::diff_forms}, refers to them as `{\em things which occur under integral signs}'. Such would not be the case if we were to use vectors, because integration of vector quantities is a metric operation. The same holds for vector operations; the grad, curl and div are metric-dependent operators, whereas the exterior derivative, which plays a similar role for differential forms, is metric-free. The important difference between vectors and forms will be explicitly addressed in this paper.

When integrals over $k$-dimensional geometric objects are considered, the orientation of these $k$-dimensional objects need to be taken into account. If we change the orientation of a point, curve, surface or volume, some integral values change sign, whereas others do not. For instance the work $W_{AB}$ of a conservative force along a curve $\mathbf{\gamma}$ connecting the points $A$ and $B$ is equal to $-W_{BA}$, i.e. the work of the same force in the opposite direction along the curve. So the physical quantity work changes sign when we change the orientation of the curve. Mass, on the other hand, which is the integral of mass density over a volume, does not change sign when we change the orientation. 

Therefore, we need to consider two distinct types of differential forms: Those that do not change sign when orientation is reversed, the {\em true forms} and those that do change sign, the {\em pseudo-forms}. The operator which switches between forms and pseudo-forms is called the {\em Hodge-$\star$ operator}. This operator depends explicitly on the metric.

Integrals and integral relations can be represented without error in terms of duality pairing between chains and cochains. The distinction between integrals of true forms and pseudo-forms requires in principle two grids: One on which we represent the integral of a true form and the other grid on which we represent the integral of a pseudo-form. The formulation obtained by employing two dual grids resembles staggered finite volume methods.

An alternative way to implement the action of the Hodge-$\star$ operator is to make use of an inner product. In this approach only one grid is required. The formulation based on a single grid approach leads to a finite element method.

In this paper we try to explain this structure in more detail and show in a specific example that the single grid approach and the dual grid approach lead to equivalent solutions. By making the clear distinction between topological concepts (integrals and discrete integral representations) and metric dependent operations (Hodge-$\star$), it is very easy to switch between orthogonal grids and curvilinear grids. It will be shown that single grid and dual grid methods give equivalent solutions in a specific problem.

Commuting relations between the discretization (mimetic projection) and operations at the continuous and discrete level will play an important role in order to ensure that the `discrete system behaves just like the continuous system'.


Throughout this paper the basic idea will be highlighted by putting statements in a box and the main idea of this paper is:

\vskip 0.3cm
\noindent
\framebox[1.0\linewidth]{
\begin{minipage}{0.92\linewidth}
A discrete representation of a physical system will display the same structure/dynamics when when discrete operators and the continuous operators commute with the projection of the infinite dimensional space onto the discrete space.
\end{minipage}}
\vskip 0.3cm

The same idea has been put forward in many different papers. Early work in this field was reported by Branin, \citep{Branin}. Dissecting physical models into metric-free components and metric-dependent parts was originally proposed by Tonti, \citep{tonti1975formal}. Application of Tonti's ideas for electromagnetism is fully treated by Mattiussi \citep{mattiussi2000finite}. A very good introduction in the geometric structure of electromagnetism is given by Bossavit, \citep{bossavit:japanese_01,bossavit1999::japan_02}. But based on the analogies described by Tonti, the construction advocated by Bossavit have a much wider range of applicability than just electromagnetism. Hodge theory of harmonic forms was described by Dodziuk, \citep{Dodziuk76}. Hyman and Scovel, \citep{HymanScovel90}, derived mimetic operators in a finite difference setting which was later generalized by Bochev and Hyman, \citep{bochev2006principles}. Mimetic finite difference methods are described in Brezzi et al. and Steinberg et al., \citep{BrezziBuffaLipnikov2009,HYmanSteinberg2004,HymanShashkovSteinberg97,RobidouxSteinberg2011,bookShashkov}. Arnold, Falk and Winther have described an extensive framework in a finite element context, \citep{arnold2006finite,arnold2010finite}. The geometric ideas underlying this paper are also extensively studied by Desbrun et al, \citep{desbrun2005discrete} and Hirani, \citep{Hirani_phd_2003} and DiCarlo et al., \citep{DiCarlo}. For fluid flow calculations, mimetic methods were used by Perot, \citep{perot2006mimetic,PerotSubramanian2007a} and the importance of preserving physical invariants was illustrated in Perot's review paper \citep{perot43discrete}. Application of these ideas for spectral elements was described in \citep{kreeft2011mimetic,Rufat} and application of these ideas to Stokes flow can be found in \citep{kreeft::stokes,Kreeft_A_priori_Stokes}. Extension to compatible isogeometric methods can be found in \citep{BuffaSangalliRivasVazquez2011} and the connection of isogeometric methods with geometry can be found in \citep{hiemstra2012high}. The relation between the Hodge matrix and mass matrices in finite element methods is discussed by Hiptmair and Tarhasaari et al., \citep{hiptmair2001discrete,tarhasaari1999some}.

The outline of this paper is as follows: In Section~\ref{Section::DifferentialGeometry} we introduce the necessary background on vector fields and differential forms. In Section~\ref{sec:Orientation} the distinction between forms and pseudo-forms is discussed. In Section~\ref{Section::AlgebraicTopology} a discrete representation of integrals is given and the analogy with the continuous forms in Section~\ref{Section::DifferentialGeometry} is presented. In Section~\ref{Section::Discretization} we describe how we can convert continuous forms in discrete forms and vice versa. In Section~\ref{sec:Poisson_volume_forms} we demonstrate how the single grid and the dual grid approach can be used to solve the Poisson equation for volume forms. This shows how the method works in practice and also shows that both solutions are equivalent. In Section~\ref{Section::Conclusions} conclusions are drawn and further applications are discussed.

\section{FORMS AND VECTORS}
\label{Section::DifferentialGeometry}

This section provides an introduction of the basic elements of differential geometry. For more detailed definitions the reader is referred to \citep{frankel,flanders::diff_forms}. We emphasize here the distinction between vectors and covectors. Just like the vector field is an extension of the vector concept to manifolds, the $1$-form is an extension of the covector to manifolds. 



As Burke \citep{burke1985applied} puts it:
\begin{quote}
\emph{Were this a mere change in notation, it would make no sense to change things. It is not a mere change in notation, however, but a basic change in the fundamental concepts. The new concepts are better for unarguable reasons: they [differential forms] correctly represent a larger symmetry group, and therefore correctly represent more features of the real world.}
\end{quote}
	

\subsection{Tangent vectors and vector fields}
Before introducing differential forms, we need to define vectors or -- more precisely -- {\em tangent vectors} in a domain $\manifold{M}$. In general, $\manifold{M}$ is a differentiable manifold. Let $\gamma(t)$ be a curve in $\manifold{M}$ parametrized by $t \in (-\epsilon,\epsilon)$, $\epsilon >0$, with $p = \gamma(0) \in \manifold{M}$. The derivative $\dot{\gamma}(0)$ is a {\em tangent vector} at the point $p \in \manifold{M}$. Note that a vector will in general not lie in the manifold $\manifold{M}$ itself. So the conventional image of a vector as an arrow connecting two points in the domain $\manifold{M}$ is inadequate. For an $n$-dimensional manifold we can define $n$ curves through the point $p \in \manifold{M}$ which produce $n$ linear independent vectors. A collection of $n$ linearly independent vectors, $\left .e_1 \right |_p,\ldots ,\left .e_n \right |_p$, at the point $p$ spans a linear vector space denoted by $\tangentspace{\manifold{M}}$. Any other vector at $p$ can then be written as a linear 
combination of these basis vectors, i.e.
\[ \left . v \right |_p = \sum_{i=1}^n a^i(p) \left . e_i \right |_p \;, \]
where the expansion coefficients $a^i(p)$ are associated to the point $p$ and the particular basis $\left. e_i \right |_p$. If we introduce a local parametrization of the point $p$, we can use the coordinate functions $x^i$ for the curves which generate a basis at $p$. Such a basis is called a {\em coordinate basis} for the tangent space $\tangentspace{\manifold{M}}$ and in this case the basis vectors are generally denoted by $\left . \partial/\partial x^i \right |_p$, or briefly $\left . \partial_i \right |_p$. For a coordinate basis a vector is represented as
\[ \left . v \right |_p = \sum_{i=1}^n b^i(p) \left . \frac{\partial}{\partial x^i} \right |_p \;.\]
If we introduce another coordinate system in which to represent $p$ locally, say $(y^1,\ldots ,y^n)$, we have
\[ \left . v \right |_p = \sum_{i=1}^n \hat{b}^i(p) \left . \frac{\partial}{\partial y^i} \right |_p \;.\]
Since any $y^i = y^i(x^1,\ldots ,x^n)$, we have
\begin{eqnarray*} 
\left . v \right |_p & = & \sum_{i=1}^n b^i(p) \left . \frac{\partial}{\partial x^i} \right |_p \\
 & = & \sum_{i=1}^n \sum_{j=1}^n b^i(p) \left .  \frac{\partial y^j}{\partial x^i}\frac{\partial}{\partial y^j} \right |_p \\
 & = & \sum_{i=j}^n \hat{b}^j(p) \left . \frac{\partial}{\partial y^j} \right |_p \;.
\end{eqnarray*}
So $\hat{b}^j = \left ( \partial y^j /\partial x^i \right ) b^i$. The fact that a change of coordinates is effectively the application of the chain rule for differentiation motivates the notation $\partial_i$ for the coordinate basis.

The construction of a vector at a point $p \in {\manifold{M}}$ can be done for all points in ${\manifold{M}}$, such that we smoothly associate with each point a vector. This construction generates vector fields. The collection of all tangent spaces is called the {\em tangent bundle}
\[ T \manifold{M} := \bigcup_{p\in \manifold{M}} \tangentspace{\manifold{M}} \;.\]

\subsection{Covectors and $1$-forms}
With any linear vector space $V$ we can associate the dual space $V^*$ of linear functionals acting on $V$, i.e.
\[ \forall \alpha \in V^*\;,\;\;\; \alpha \,:\, V \rightarrow \mathbb{R}\;.\]
So with $\tangentspace{\manifold{M}}$, we can associate the dual space $T_p^*\manifold{M}$ of linear functionals acting on vector $v \in \tangentspace{\manifold{M}}$.
%
\[
\alpha: \tangentspace{\manifold{M}}\mapsto \mathbb{R} \;.
\]
\[
\alpha\left(a{v} + b{u}\right) = a\alpha({v}) + b\alpha({u})\;,\quad \forall u,v \in T_p\manifold{M} \mbox{ and }\alpha \in T_p^*\manifold{M} \;.
\]
The dual space itself becomes a linear vector space if we set
\[
\left(a\alpha + b\beta\right)\left({v}\right) = a\alpha\left({v}\right) + b\beta\left({v}\right) \;,\quad \forall u, v \in T_p\manifold{M} \mbox{ and }\alpha, \beta \in T_p^*\manifold{M} \;.
\]
The elements $\alpha,\beta \in T_p^*\manifold{M}$ are called {\em covectors at the point $p \in \manifold{M}$}. The dimension of the cotangent space $\mbox{dim} T_p^*\manifold{M} = \mbox{dim} T_p\manifold{M} = n$.

Let $\left .e_1 \right |_p,\ldots ,\left .e_n \right |_p$ be a basis for the tangent space $\tangentspace{\manifold{M}}$, then a canonical basis for the cotangent space $T_p^*\manifold{M}$ is given by $\left .e^1 \right |_p,\ldots ,\left .e^n \right |_p$, where the basis covectors satisfy
\[ e^i (e_j) = \delta_j^i = \left \{ \begin{array}{ll}
1 \quad & \mbox{if } i=j \\
 & \\
0 & \mbox{if } i \neq j
\end{array} \right .
\;.\]

In case the vectors are represented in a {\em coordinate basis}, $\partial_i$, the basis covectors are denoted by $dx^j$. In this case any covector at the point $p$ can be represented as
\[ \alpha = \alpha_1(p) \left . \ederiv x^1 \right |_p + \ldots + \alpha_n(p) \left .\ederiv x^n \right |_p \;.\]
Application of such a covector to a tangent vector at $p$ then yields
\begin{eqnarray}
\alpha (v) & = & \sum_{i=1}^n \alpha_i(p) \ederiv x^i \left ( \sum_{j=1}^n v^j(p) \frac{\partial}{\partial x^j} \right ) \nonumber \\
 & = & \sum_{i=1}^n \sum_{j=1}^n \alpha_i(p) v^j(p) \ederiv x^j \left (\frac{\partial}{\partial x^j} \right ) \nonumber \\
 & = & \sum_{i=1}^n \alpha_i(p) v^i(p) \;. \label{eq:duality_pairing_vector_covector}
\end{eqnarray}

\begin{remark}
In a given basis, application of a covector to a vector resembles the inner product in Euclidean space. However, an inner product depends on the metric tensor, whereas $\alpha(v)$ is independent of the metric. Furthermore, an inner product is a bilinear form on a single vector space $V$, whereas application of a covector to a vector is an operation between two distinct spaces: $V^* \times V \rightarrow \mathbb{R}$. That these spaces are really different can be seen when we apply a change of coordinates. In this case vectors and covectors transform differently and should therefore be treated as different mathematical entities. 
\end{remark}

\begin{remark}
A pictorial representation of a vector is usually in terms of an arrow. Covectors can be represented by $(n-1)$-dimensional hyper-surfaces, (ordinary surfaces in 3D). Duality pairing between a covector and a vector then yields the number of times the arrow (vector) pierces the surfaces (covectors). For examples of these graphical representations, see Bossavit, \citep{bossavit:japanese_01}, Burke, \citep{burke1985applied} and Misner, Thorne and Wheeler, \citep{misner1973gravitation}. The ultimate aim of this distinctive representation is to emphasize the difference between forms and vectors.
\end{remark}

\begin{remark}
The cotangent space $T_p^*\manifold{M}$ is isomorphic to the tangent space $T_p\manifold{M}$. But there is no canonical isomorphism which associates elements $\alpha \in T_p^*\manifold{M}$ to elements $v \in T_p\manifold{M}$. It is {\em physics} which provides a unique connection between vectors and covectors.
\end{remark}

\begin{remark}
We cannot 'see' a form, i.e. a linear functional, directly, but we can only assess its action on the elements of the primal space. The only knowledge we can obtain of a form $\alpha$ is by applying $\alpha$ to various elements $v$ of the vector space. A similar thing happens for some physical variables: Nobody has seen 'force', but we only see the action of a force. We use scales and spring balances to measure force. The deflection of the hands on the scale or the extension of the spring are used to measure force. So we only have access to the work performed by the force, $F(v)$. It therefore seems obvious to represent certain physical variables by forms.
\end{remark}

Just as we extended the tangent space in a point $p \in \manifold{M}$, to all points in $\manifold{M}$ to form the tangent bundle, we can also consider the collection of cotangent spaces for all points in $\manifold{M}$. This defines the cotangent bundle
\[ T^* \manifold{M} := \bigcup_{p\in \manifold{M}} T_p^* \manifold{M} \;.\]
An element (section) of $T^* \manifold{M}$ is then written as
\[ \kdifform{\alpha}{1} = \sum_{k=1}^{n}\alpha_{k}(x^1,\ldots , x^n) \ederiv x^{k} \;.\]
Such an element from the cotangent bundle is called a {\em differentiable $1$-form} or a {\em $1$-form}. The space of differentiable $1$-forms is also denoted by $\Lambda^1(\manifold{M})$.

Application of a $1$-form $\kdifform{\alpha}{1}$ to a vector field $v$ assigns to every point in $\manifold{M}$ a real number, i.e. $\kdifform{\alpha}{1} ( v)$ is a real-valued function on $\manifold{M}$. We will denote the space of real-valued functions on $\manifold{M}$ by $\Lambda^0(\manifold{M})$, the space of $0$-forms on $\manifold{M}$.

%

\begin{example}\label{ex:line_integral}
Consider a 1-form $\kdifform{\varphi}{1}(x^{1},x^{2}) = \varphi_{1}(x^{1},x^{2})\ederiv x^{1} + \varphi_{2}(x^{1},x^{2})\ederiv x^{2}$ and a curve $\boldsymbol{\gamma}(s) = \left(\gamma^{1}(s),\gamma^{2}(s)\right)$, such that:
           \[
              \boldsymbol{\gamma}:[0,1]\mapsto\manifold{M}
           \]
           The tangent vectors along the curve are given, as usual, by:
           \[
              {\vec{g}}(s) = \frac{\ederiv \gamma^{1}}{\ederiv s} \,\partial_{1} + \frac{\ederiv \gamma^{2}}{\ederiv s}\,\partial_{2}
           \]
           The action of $\kdifform{\varphi}{1}$ on ${\vec{g}}$ is then given by:
           \[
              \kdifform{\varphi}{1}\left({\vec{g}}\right) = \sum_{k=1}^{n}\varphi_{k}\frac{\ederiv\gamma^{k}}{\ederiv s} = w(s)
           \]
           Where $w(s)$ is a function of $s$ and whose Lebesgue integral is:
           \[
              \int_{0}^{1}w(s) = W
           \]
\end{example}

\begin{remark}           
Note that the less usual notation for Lebesgue integral was used. Typically the notation for Lebesgue integral shows the measure term, in this case $\ederiv s$, but that is optional. The Lebesgue measure was omitted since, by historical misfortune, $\ederiv s$, the Lebesgue measure has the same notation as a 1-form basis $\ederiv x^{1}$, nevertheless both are distinct mathematical objects. This overload of notation is the root of a common confusion between differential and differential forms (a functional) and measures: A form changes sign under a change of orientation, whereas a measure remains positive.

Many integrals in physics are still considered as Riemann integrals, where the $\ederiv s$ is interpreted as the limit $\Delta s \rightarrow \ederiv s$. But the Riemann integral is defined only for a limited class of functions and this set of functions is not closed under the operation of taking point-wise limits of sequences of functions in this class, \citep{LiebLoss}. All integrals will be Lebesgue integrals in this paper.
\end{remark}

\begin{remark}
As for the physical meaning of this example, consider $\kdifform{\varphi}{1}$ as a force. Hence, $W$ is the work done by this force along the path $\boldsymbol{\gamma}(s)$ and $w(s)$ is a density of work along the path. The integral along a curve is a metric-free operation when the force is represented as a differential form. If the force were represented as a vector, as is the case in many elementary physics books, work would be a metric-dependent concept. 
\end{remark}

\begin{remark}
Duality between vectors and covectors is reflexive and therefore we can interpret Example~\ref{ex:line_integral} as a recipe to integrate the $1$-form $\kdifform{\varphi}{1}$ along the path $\boldsymbol{\gamma}$, or to integrate the path $\boldsymbol{\gamma}$ against the $1$-form $\kdifform{\varphi}{1}$, \citep{tao2007differential}.
\end{remark}

Now that we have defined the $0$-forms on $\manifold{M}$ and the $1$-forms on $\manifold{M}$, we can proceed in different ways to define $k$-forms for $k>1$. We can either define exterior powers of a vector space to produce $k$-vectors and then define $k$-forms as the elements of the dual space of the space of $k$-vectors. This approach is described in \citep{tonti1975formal, flanders::diff_forms, rwr1994differential}. Alternatively, we can define a $k$-form as an alternating $k$-tensor on the tangent space which maps into the real numbers, \citep{spivak1998calculus}
\[ \kdifform{\alpha}{k}\,:\, \underbrace{T_p\manifold{M} \times \dots \times T_p\manifold{M}}_{k \mbox{ copies}} \rightarrow \mathbb{R} \;,\]
with 
\[ \kdifform{\alpha}{k} (\ldots, v_i ,\ldots,v_j,\ldots ) = - \kdifform{\alpha}{k} (\ldots, v_j ,\ldots,v_i,\ldots ) \;,\;\;\; v_i \in T_p\manifold{M}\;,\;\; i=1,\ldots,k\;.\]

Here we can define the {\em exterior product} or {\em wedge product} to inductively construct $k$-forms, for $1 \leq k \leq n$. Let $\Lambda^k(\manifold{M})$ and $\Lambda^l(\manifold{M})$ the space of $k$-forms and $l$-forms, respectively, with $k+l\leq n$, then the \emph{wedge product}, $\wedge$, is a mapping:
\[
	\wedge : \differentialFormSpace{k}{\manifold{M}} \times \differentialFormSpace{l}{\manifold{M}} \rightarrow  \differentialFormSpace{k+l}{\manifold{M}} , \quad k+l \le n
\]
that satisfies the following properties:
\begin{subequations}
			\begin{align}
				(\kdifform{\alpha}{k} + \kdifform{\beta}{l})\wedge\kdifform{\gamma}{m} = \kdifform{\alpha}{k}\wedge\kdifform{\gamma}{m} + \kdifform{\beta}{l}\wedge\kdifform{\gamma}{m} & \quad \text{(Distributivity)} \\
				(\kdifform{\alpha}{k}\wedge\kdifform{\beta}{l})\wedge\kdifform{\gamma}{m} = \kdifform{\alpha}{k}\wedge(\kdifform{\beta}{l}\wedge\kdifform{\gamma}{m}) & \nonumber \\
				 = \kdifform{\alpha}{k}\wedge\kdifform{\beta}{l}\wedge\kdifform{\gamma}{m} & \quad \text{(Associativity)} \label{wedgeassociativity}\\
				a \kdifform{\alpha}{k}\wedge\kdifform{\beta}{l} = \kdifform{\alpha}{k}\wedge a\kdifform{\beta}{l} = a (\kdifform{\alpha}{k}\wedge\kdifform{\beta}{l}) & \quad \text{(Multiplication by functions)} \\
				\kdifform{\alpha}{k}\wedge\kdifform{\beta}{l} = (-1)^{kl}\kdifform{\beta}{l}\wedge\kdifform{\alpha}{k}& \quad \text{(Skew symmetry)} \label{wedge::skew_symmetry}
			\end{align}
\end{subequations}
where $\kdifform{\alpha}{k}\in\Lambda^{k}(\manifold{M})$, $\kdifform{\beta}{l}\in\Lambda^{l}(\manifold{M})$ and $\kdifform{\gamma}{m}\in\Lambda^{m}(\manifold{M})$.

Consider a sufficiently smooth bounded $n$-dimensional oriented manifold $\manifold{M} \subset \realNumber^{n}$ with boundary $\partial \manifold{M}$ and a coordinate system where points in the manifold, $\mathbf{x}\in \manifold{M}$, are represented by an $n$ tuple $\mathbf{x}:= \left(x^{1}, \cdots,  x^{n}\right)$. Let $\alpha^{(k)}$ denote a differential $k$-form, $k\leq n$ and  $\differentialFormSpace{k}{\manifold{M}}$ denote the space of \emph{differential k-forms} or \emph{k-forms}, then $\alpha^{(k)} \in \differentialFormSpace{k}{\manifold{M}}$, can be written as

\begin{align}
\alpha^{(k)}= \sum_{I} \alpha_{I}\left( \mathbf{x} \right) dx^{i_{1}} \wedge dx^{i_{2}} \wedge \cdots \wedge dx^{i_{k}} ,
\label{equation::differentialForms_definition}
\end{align}
where $I = i_{1}, \cdots, i_{k},$ and $1 \leq i_{1} < \cdots < i_{k} \leq n$ and where $\alpha_{I}\left( \mathbf{x} \right)$ are continuously differentiable scalar functions. 


\subsection{Integration of differential forms}
\vskip 0.3cm
\noindent
\framebox[1.0\linewidth]{
\begin{minipage}{0.92\linewidth}
A Physical system generally preserves integral quantities and these quantities are topological in the sense that they do not depend on a particular coordinate system or the metric employed. The integration of differential forms is independent of any metric notions in space, \citep[Ch.3]{frankel}; no Riemann metric, dot products or $k$-dimensional volume elements are required. If we would have used vectors for integration all these metric concepts need to be taken into account. Integrals are of paramount importance in mimetic discretizations, because we cannot represent differential forms and the operations on and between differential forms in a finite dimensional setting. But we {\em can} represent integrals of differential forms and the integral relations between them exactly.
\end{minipage}}
\vskip 0.3cm

In $\realNumber^{3}$ there are four types of differential forms, that is 0-,1-, 2- and 3-forms, as many as the different kinds of submanifolds: i.e. points, lines, surfaces and volumes. Differential $k$-forms naturally integrate over $k$-dimensional manifolds. Let $\alpha^{(k)} \in \differentialFormSpace{k}{\manifold{M}}$ and $\manifold{M}_k \subset \realNumber^{n}$, with $k=dim(\manifold{M}_{k})$,
\begin{align}
\duality{\alpha^{(k)}}{\manifold{M}_{k}} :=\int_{\manifold{M}_{k}} \alpha^{(k)},
\label{eq:integration_difform}
\end{align}
this represents a metric-free duality pairing.

Differential forms live on manifolds and can be easily transformed from one manifold to another, under the action of a special mapping called the {\em pullback}. Let $\Phi : \manifold{M}_{ref} \rightarrow \manifold{M}$ be a mapping between two manifolds, the pullback operator, $\pullback : \differentialFormSpace{k}{\manifold{M}} \rightarrow \differentialFormSpace{k}{\manifold{M}_{ref}}$ is such that,
\begin{align}
\int_{\Phi\left(\manifold{M}_{ref} \right)} \alpha^{(k)} = \int_{\manifold{M}_{ref}} \pullback \alpha^{(k)}.
\end{align} 

The actual evaluation of the integrals is analogous to the case discussed in Example~\ref{ex:line_integral}. Choose a parameterization of the manifold $\manifold{M}_k$
\[ \tau_k \,:\, [0,1]^k \rightarrow \manifold{M}_k \;,\]
and then
\[ \int_{\manifold{M}_{k}} \alpha^{(k)} = \int_{\tau_k([0,1]^k)} \alpha^{(k)} = \int_{[0,1]^k} \tau_k^\star (\alpha^{(k)}) \;,\]
where the integral on the right denotes the usual integral in $\mathbb{R}^k$. For more complex manifolds one might use a collection of parameterizations of the manifold, $\tau_k^i$, for the representation of the manifold in which case the integral is evaluated as
\begin{equation}
\int_{\manifold{M}_{k}} \alpha^{(k)} = \sum_i \int_{\tau_k^i([0,1]^k)} \alpha^{(k)} = \sum_i \int_{[0,1]^k} (\tau_k^i)^\star (\alpha^{(k)}) \;.
\end{equation}



The inner-product of $k$-forms, $\inner{\cdot}{\cdot}$, is a symmetric, bi-linear mapping:
\[
	\inner{\cdot}{\cdot}: \differentialFormSpace{k}{\manifold{M}} \times \differentialFormSpace{k}{\manifold{M}} \rightarrow \Lambda^{0}(\manifold{M})\;.
\]
Let $\kdifform{\alpha}{k}$ and $\kdifform{\beta}{k}$ be two $k$ forms written in the form (\ref{equation::differentialForms_definition}), then the point-wise inner product is given by, \citep{frankel,flanders::diff_forms,rwr1994differential}
\begin{equation}
\left ( \kdifform{\alpha}{k},\kdifform{\beta}{k} \right ) = g^{i_1.j_1}\dots g^{i_k,j_k} \alpha_{i_1,\dots,i_k} \beta_{j_1,\dots,j_k} \;,
\end{equation}
where $g^{k,l}:= (\ederiv x^k,\ederiv x^l)$ are the components of the metric tensor.

A combination between the wedge operator and the inner product induces the Hodge-$\star$ operator, $\star$, which is a mapping:
\begin{equation}
	\star: \differentialFormSpace{k}{\manifold{M}} \rightarrow \differentialFormSpace{n-k}{\manifold{M}}
	\label{eq:SpacesHodge}
\end{equation}
such that:
\begin{align}
	\alpha^{(k)} \wedge \star \beta^{(l)} := \inner{\alpha^{(k)}}{\beta^{(k)}} \sigma^{(n)},
\label{eq:HodgeInnerProduct}	
\end{align}
where $\sigma^{(n)} \in \differentialFormSpace{n}{\manifold{M}}$ is the volume form, that is the $n$-differential form defined on an $n$-manifold $\manifold{M}$ such that $\int_{\manifold{M}}\kdifform{\sigma}{n}$ is the volume of the manifold $\manifold{M}$ and we have $\star 1 = \kdifform{\sigma}{n}$. The important role of the Hodge-$\star$ with regard to orientation will be explained in the next section.

The $L^{2}$ inner product, $\innerspace{\cdot}{\cdot}{L^{2}(\manifold{M})}$, is defined as a mapping:
\[
	\innerspace{\cdot}{\cdot}{L^{2}(\manifold{M})}:\differentialFormSpace{k}{\manifold{M}} \times \differentialFormSpace{k}{\manifold{M}} \rightarrow \realNumber\;, 
\]
such that:
\begin{align}
	\innerspace{\alpha^{(k)}}{\beta^{(k)}}{L^{2}(\manifold{M})} := \int_{\manifold{M}} \left( \alpha^{(k)}, \beta^{(k)} \right)\kdifform{\sigma}{n} = \int_{\manifold{M}} \alpha^{(k)} \wedge \star \beta^{(l)}\;.
\label{eq::diffGeom_L2_inner}
\end{align}

A crucial operator in differential geometry is the exterior derivative, $\ederiv$, defined as a mapping
\[ 
	\ederiv: \differentialFormSpace{k}{\manifold{M}} \rightarrow \differentialFormSpace{k+1}{\manifold{M}}
\] 
\vskip 0.3cm
\noindent
\framebox[1.0\linewidth]{
\begin{minipage}{0.92\linewidth}
The next equation is the most important equation in mimetic methods and can therefore be referred to as {\em the Mother of all equations}. It relates the geometric boundary operator to the exterior derivative.
\[  \mbox{The generalized Stokes Theorem} \]
\begin{align}
	\int_{\manifold{M}_{k+1}} d \alpha^{(k)}= \int_{\partial \manifold{M}_{k+1}} \alpha^{(k)}
	\quad \Leftrightarrow \quad \left<d \alpha^{(k)},\manifold{M}_{k+1}  \right> = \left< \alpha^{(k)}, \partial \manifold{M}_{k+1} \right>,
	\label{eq:generalized_Stokes}
\end{align} 
\end{minipage}}
\vskip 0.3cm

where $\manifold{M}_{k+1}$ is a $(k+1)$-dimensional manifold and $\partial\manifold{M}_{k+1}$ is its boundary, a $k$-dimensional manifold, and $\kdifform{\alpha}{k} \in \differentialFormSpace{k}{\manifold{M}}$. This represents the  \emph{generalized Stokes theorem} that encapsulates the Newton-Leibniz, Stokes and Gauss theorem and hence generalizes the vector calculus $\nabla$, $\nabla\times$ and $\nabla\cdot$ operators to arbitrary manifolds. Moreover, this operator is independent of any metric or coordinate system and satisfies the Leibniz rule $\ederiv \left( \alpha^{(k)} \wedge \beta^{(l)} \right) = \ederiv \alpha^{(k)} \wedge \beta^{(l)} + (-1)^{k} \alpha^{(k)} \wedge d \beta^{(l)}$, and is nilpotent, $\ederiv\ederiv\alpha^{(k)} := 0$. The pullback operator and the exterior derivative have the following commuting property,
\begin{equation*}
\begin{tikzpicture}[descr/.style={fill=white,inner sep=2.5pt}]
\matrix (m) [matrix of math nodes, row sep=1.5em,
column sep=1.5em]
{\differentialFormSpace{k}{\manifold{M}} & \differentialFormSpace{k+1}{\manifold{M}} \\
\differentialFormSpace{k}{\manifold{M}_{ref}} & \differentialFormSpace{k+1}{\manifold{M}_{ref}}  \\};
\path[->,font=\scriptsize]
(m-1-1) edge node[auto] {$d$} (m-1-2)
(m-2-1) edge node[auto] {$d$} (m-2-2)
(m-1-1) edge node[auto] {$\pullback$} (m-2-1)
(m-1-2) edge node[auto] {$\pullback$} (m-2-2);
\end{tikzpicture}
\end{equation*}

Finally, the last piece of machinery is the co-differential, $\coderiv$, which is a mapping:
\[
	\coderiv: \differentialFormSpace{k}{\manifold{M}} \rightarrow \differentialFormSpace{k-1}{\manifold{M}}\;,
\]
that satisfies the identity:
\begin{align}
\left( \alpha^{(k-1)}, \coderiv \beta^{(k)} \right)_{\manifold{M}} = \left(d \alpha^{(k-1)}, \beta^{(k)} \right)_{\manifold{M}} - \int_{\partial \manifold{M}} \mbox{tr}\, \alpha^{(k-1)} \wedge \mbox{tr}\, \star \beta^{(k)}\;,
\label{eq:integration_by_parts}
\end{align}
where tr is the trace operator, which is the pullback of the immersion of boundary of the manifold into the manifold, $\iota^*$, with $\iota\,:\,\partial \manifold{M} \mapsto \manifold{M}$.
When $\int_{\partial \manifold{M}} \mbox{tr}\, \alpha^{(k-1)} \wedge \mbox{tr}\, \star \beta^{(k)}=0$ then $\coderiv$ is the formal Hilbert adjoint of the exterior derivative induced by the $L^{2}$ inner product  \eqref{eq::diffGeom_L2_inner}:
\begin{align}
\innerspace{d \alpha^{(k-1)}}{\beta^{(k)}}{\manifold{M}} = \innerspace{\alpha^{(k-1)}}{\coderiv \beta^{(k)}}{\manifold{M}},
\end{align} 
Opposite to the exterior derivative, the codifferential is a metric-dependent operator. given by, $\coderiv = \left(-1 \right)^{n(k+1)+1}\star d \star$.

The Laplace operator, $\Delta$, is a mapping:
\[
	\Delta: \differentialFormSpace{k}{\manifold{M}} \rightarrow \differentialFormSpace{k}{\manifold{M}}
\]
which, in terms of the exterior derivative and the co-differential, is given by
\begin{align}
 \Delta = \coderiv \ederiv + \ederiv \coderiv .
 \label{eq:Laplace-DeRham}
\end{align}

For a $0$-form $\kdifform{\alpha}{0} \in \differentialFormSpace{0}{\manifold{M}}$ the Laplace operator is given by $\Delta \kdifform{\alpha}{0} = \coderiv \ederiv \kdifform{\alpha}{0}$. For $n$-forms $\kdifform{\omega}{n} \in \differentialFormSpace{n}{\manifold{M}}$ the Laplace operator is $\Delta \kdifform{\omega}{n} = \ederiv \coderiv \kdifform{\omega}{n}$.

\section{ORIENTATION}
\label{sec:Orientation}
Now that we have introduced sub-manifolds and differential forms as metric-free ingredients in the evaluation of integrals, we need to discuss orientation of space and its geometric objects. Orientation is something we need to introduce to do our calculations, but physics should -- in a sense -- be independent of the choice of our orientation. 

From our earliest physics courses orientation has played a role: the right-hand rule in electromagnetism, the outward unit normal, calculation of circulation counter-clockwise, vorticity in aerodynamics is oriented clockwise, definition of the vector product of two vectors, etc. Many sign conventions in physics are also due to a preferred orientation. 

The basic question with orientation is: `What happens in the description of physics if we change the orientation?' Of course, it will only affect the description of physics and not physics itself. A few examples might illustrate the point.

\begin{example}
In Example~\ref{ex:line_integral} we calculated the line integral along the path $\mathbf{\gamma}(s) = (\gamma^1(s),\gamma^2(s))$, $0\leq s \leq 1$. If the $1$-form integrated along this curve denotes force, then the line integral denotes the amount of work, $W$, done by the force along the curve.

Now suppose that we change the direction (orientation) in which we traverse the curve, i.e. we now take the curve $\mathbf{\gamma}(s) = (\gamma^1(1-s),\gamma^2(1-s))$, $0\leq s \leq 1$. Then the value of this integral is $-W$ for a conservative force. This is as expected. So changing the orientation in the calculation of work for a conservative force leaves the differential form unchanged and reverses the integral: $W_{AB}=-W_{BA}$, if the point $A$ and $B$ are the endpoints of the curve on the manifold $\manifold{M}$.
\end{example}

Is this generally true? If we change the orientation of the domain of integration, does the integral quantity always change sign? The answer is no. A simple counter example is the $3$-form mass density, $\kdifform{\rho}{3}$. 

If we integrate $\kdifform{\rho}{3}$ over a positively oriented volume (right-hand rule), we obtain the mass in that volume
\[ M = \int_{\manifold{M}_3^+} \kdifform{\rho}{3} \geq 0 \;.\]
If we change the orientation of the volume and the differential form $\kdifform{\rho}{3}$ would remain unchanged, the integral would be
\[ M = \int_{\manifold{M}_3^-} \kdifform{\rho}{3} = - \int_{\manifold{M}_3^+} \kdifform{\rho}{3} \leq 0 \;.\]
So a change in orientation would lead to negative mass, which is physically unacceptable. Not because we think of mass as being positive -- hypothetically mass can be negative --, but because its sign depends on the orientation. In order to restore our faith in physics and disappoint those who suffer from overweight, the mass density $3$-form should change sign when the orientation is changed. In this case, mass will always be positive, regardless of the sense of orientation.

So we have (at least) two classes of differential forms: those that do not change sign when the orientation is reversed and those that do change sign when the sense orientation is changed.

\subsection{Orienting geometric objects}
Orientation is an indispensable part the description of physics and physics compatible numerical techniques should take orientation into account; either explicitly by assigning an orientation to the computational mesh, or implicitly by following sign conventions. Orientation is discussed in any textbook on differential geometry such as \citep[\S 2.8]{frankel}. A very good introduction is also given by Bossavit, \citep{bossavit1999::japan_02}.

As we will show in the remainder of the paper, orientation is also an important aspect of computational physics.

\subsubsection{Orientation of a vector space}
Let $V$ be a vector space and let $\mathbf{e} = (\mathbf{e}_1,\ldots ,\mathbf{e}_n)$ and $\mathbf{f} = (\mathbf{f}_1,\ldots ,\mathbf{f}_n)$ be two bases. Then there exists a non-singular transition matrix $P$ such that $\mathbf{e} = P \mathbf{f}$. If det$(P)$ is positive, then both bases have the same orientation and if det$(P)$ is negative both bases have opposite orientation. This divides all bases into two equivalence classes. All bases in a class have the same orientation. This equivalence relation is reflexive, transitive and symmetric.

We orient a vector space by selecting a basis from one of the two equivalence classes and call this the {\em positive} orientation or the {\em default} orientation. There is nothing intrinsically 'positive' about this basis, it is just a choice. If we take any basis from the same equivalence class, this basis will have the same orientation. If we take a basis from the other class it will have the opposite orientation. Denote the equivalence class from which we pick out default orientation $Or$ and the other class $-Or$. Bossavit calls the orientation from the default orientation {\em direct frames} and the bases from $-Or$ {\em skew frames}.

\subsubsection{Orientation of a geometric object}
Consider the basis geometric objects in 3D, the point, curve, surface and volume. We can consider these objects as sub manifolds in 3D and therefore they possess charts from the objects to $\mathbb{R}^k$, $k=0,1,2,3$. Then we set up a coordinate basis in each chart. All coordinate systems within a chart have the same orientation. At the overlap between charts, we impose that the Jacobian is positive thus ensuring that in a neighboring patch the orientation is in the same equivalence class. If we can do this consistently for the whole manifold, we call the manifold {\em orientable}. Not all manifolds are orientable. A notable example is the M\"{o}bius strip and the projective plane $\mathbb{R}P^2$.

\subsection{Forms and pseudo-forms}
In the examples above, work and mass, we saw that some variables change sign when the orientation is reversed, while others do not change sign. The forms which do not change sign are called {\em genuine forms}, while forms which do change sign are called {\em pseudo-forms}. This is an important distinction between physical variables. We can never equate a form to a pseudo form, because a change of orientation would falsify the relation. Before we can equate a form to a pseudo form, either the form needs to be 'translated' to a pseudo-form, or the pseudo-form needs to be converted to a form.

The notions form and pseudo-form are taken from differential geometry and can be found in almost any book on the subject. In this paper we loosely follow Frankel, \citep[Ch.2]{frankel} and Bossavit, \citep{bossavit:japanese_01}. One may find other classifications in literature which mean exactly the same thing: Tonti, \citep{tonti1975formal}, describes the two classes as {\em configuration variables} and {\em source variables}. In earlier work, \citep{kreeft2011mimetic} we referred to {\em inner-oriented variables} and {\em outer-oriented variables}. The genuine forms were referred to as inner-oriented variables and the pseudo-forms were called outer-oriented. Bossavit, \citep{bossavit:japanese_01} talks about {\em axial vectors} or {\em twisted vectors} and {\em true vectors}. Whatever the names one uses to refer to the two distinct types of forms, it is important to acknowledge this difference. Not only for the sake of a sound physical description, but also for a proper numerical model.

\vskip 0.3cm
\noindent
\framebox[1.0\linewidth]{
\begin{minipage}{0.92\linewidth}
The operator which converts forms into pseudo-forms and vice versa is the Hodge-$\star$ operator, \citep[\S 14.1]{frankel}, defined in (\ref{eq:SpacesHodge}) and (\ref{eq:HodgeInnerProduct}).
\end{minipage}}
\vskip 0.3cm

\section{DISCRETE REPRESENTATION OF INTGRALS}
\label{Section::AlgebraicTopology}
\begin{figure}[ht!]
\centering
\includegraphics[width=0.85\textwidth]{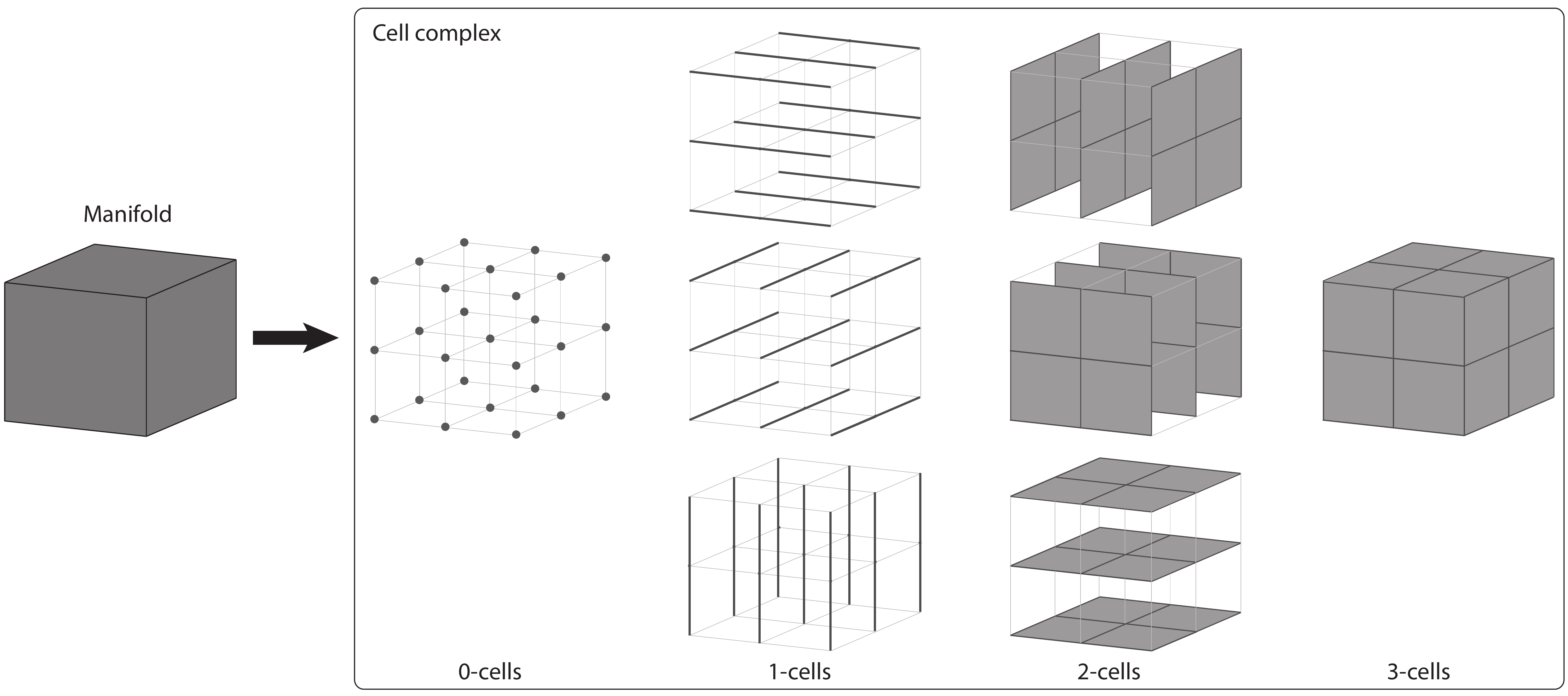}
\caption{Three dimensional domain (left) partitioned into a collection of points, lines segments, surfaces and volumes (right).}
\label{fig:cell_complex_final}
\end{figure}

Consider a three dimensional domain $\manifold{M}$ and its associated grid, as shown in Figure~\ref{fig:cell_complex_final}.
The grid not only consists of points as is common in many numerical methods, but also the line segments connecting the points, surfaces bounded by these line segments and volumes bounded by these surfaces. The partitioning of the domain (manifold) in these geometric building blocks is an instantiation of a cell complex. Loosely speaking a cell complex consists of a collection of sets, $C_k$, in this case the $k$-dimensional objects in the grid, and an operator $\partial$ which maps elements of the set $C_k$ into the set $C_{k-1}$, such that $\partial \partial C_k = \emptyset$. For a grid to be considered as a cell complex we impose that if a $k$-dimensional geometric object is in the grid, its boundary is also in the grid, where the operator $\partial$ is the boundary operator.

If we endow all geometric objects (sub-manifolds) with a default orientation, then we call the cell complex an {\em oriented cell complex}.

Let $\manifold{M}$ be a manifold covered by an oriented cell complex, $D$, or an \emph{oriented grid} consisting of points, lines, surfaces and volumes. We will call the individual $k$-dimensional geometric objects in $D$ {\em $k$-cells} $\tau_{(k),i}$, $k = 0, \cdots, n$, where $k$ denotes the dimension of the object and $i$ is a label to distinguish the different points and lines, etc. Given the cell complex $D$, the space of $k$-chains of D, $C_{k}\left(D\right)$, is the collection of weighted $k$-cells. A $k$-chain, $\kchain{c}{k} \in C_{k} \left( D \right)$ is a formal sum of $k$-cells, $\tau_{(k),i} \in D$,
\begin{align}
\kchain{c}{k} = \sum_{i} c^{i} \tau_{(k),i}.
\end{align}
By formal sum we mean a collection of cells and weight $\{ \tau_{(k),i},c^i\}$, where the $c^i$ denote the weights. The addition of two geometric objects, say points, is not defined, therefore this weighted collection, the formal sum, should not be confused with ordinary summation.

The boundary operator on $k$-chains, $\partial : C_{k} \left(D \right) \rightarrow  C_{k-1} \left(D \right)$, is an homomorphism,
\begin{align}
\partial \kchain{c}{k} = \partial \sum_{i} c^{i} \tau_{(k),i} := \sum_{i} c^{i} \partial \left(  \tau_{(k),i} \right).
\end{align}

The boundary of a $k$-cell, $\tau_{(k),i}$ is a $(k-1)$-chain formed by the oriented faces of $\tau_{(k),i}$. 
The coefficients of this $(k-1)$-chain associated to each of the faces is given by the orientations. 

\begin{align}
 \partial  \tau_{(k),i} = \sum_{j} e_{i}^{j} \tau_{(k-1),j},
\end{align}

with
\begin{align*}
\left\{
 \begin{array}{c l}
  e_{i}^{j} = 1 & \text{ , if the orientation of $\tau_{(k-1),j}$ equals the default orientation} \\
  e_{i}^{j} = -1 & \text{ , if the orientation of $\tau_{(k-1),j}$ is opposite the default orientation} \\
  e_{i}^{j} = 0 & \text{ , if $\tau_{(k-1),j}$ is not a face of $\tau_{(k),j}$} 
 \end{array}
\right.
\end{align*}

Once basis $k$-cells have been chosen, the space of $k$-chains, $\kchainspacedomain{k}{\ccomplex{D}}$, can be represented by a column vector containing only the coefficients, $c^i$, of the chain. That is, there is an isomorphism $\psi$:
\begin{equation}
\psi: \kchainspacedomain{k}{\ccomplex{D}}\mapsto\mathbb{R}^{p}, \quad p=\mathrm{rank}(\kchainspacedomain{k}{\ccomplex{D}})\;,\label{eq::algTop_isomorphism_map}
\end{equation}
defined by
\begin{equation}
\psi(\kchain{c}{k}) = \psi \left( \sum_{i}c^{i}\tau_{(k),i} \right) =  [c^{1} \cdots c^{p}]^{T}, \quad p=\text{rank}(\kchainspacedomain{k}{\ccomplex{D}})\;,
\label{eq::algTop_isomorphism_definition}
\end{equation}
where the rank of $\kchainspacedomain{k}{\ccomplex{D}}$ is the number of $k$-cells in the cell complex $D$ and the $c^i$ are the coefficients of the $k$-chain $\mathbf{c}_{(k)}$. The $k$-chain, $\kchain{c}{k}$, is printed in boldface, whereas the vector $c_{(k)}$ of coefficients is printed in regular face.

\begin{remark}
Instead of taking the individual $k$-cells, $\tau_{(k),i}$, as basis, one may choose to take any other collection of linearly independent $k$-chains, $\tilde{c}_{(k),j}$, as basis. Both bases are related by $\tau_{(k),i}=A_i^j \tilde{c}_{(k),}$, we have
\[ c_{(k)} = \sum_i c^i \tau_{(k),i} = \sum_i c^i \sum_j A_i^j \tilde{c}_{(k),j} = \sum_j \tilde{c}^j \tilde{c}_{(k),j} \;,\]
where $\tilde{c}^j = A_i^j c^i$. This change of basis allows one to incorporate isogeometric methods in this mimetic framework, \citep{hiemstra2012high}.
\end{remark}

\begin{example}\label{ex:incidence_matrix}
Consider as an example the 2-dimensional grid shown in Figure~\ref{fig:Simple_2D_grid}. The points $P_i$ are assumed to be sinks. So when something 'enters' a point we  denote this by the positive orientation and when something 'leaves' the point it will have the opposite orientation. The default orientation for the lines and surface is indicated in the figure.
\begin{figure}[ht!]
\centering
\includegraphics[width=0.35\textwidth]{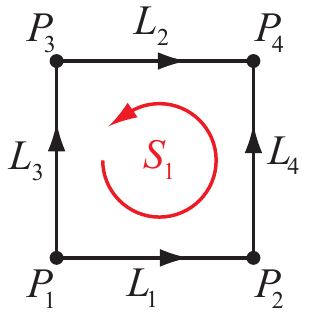}
\caption{Small oriented 2-dimensional cell complex}
\label{fig:Simple_2D_grid}
\end{figure}
For the figure we see that the boundary of the line $L_1$ is given by $\partial L_1 = P_2 - P_1$. Similarly, we can write down the boundaries for all line segments and collect everything in a matrix
\begin{equation}
\left ( \begin{array}{c} 
\partial L_1 \\
\partial L_2 \\
\partial L_3 \\
\partial L_4
\end{array} \right ) = \left ( \begin{array}{cccc}
-1 & 1 & 0 & 0 \\
0 & 0 & -1 & 1 \\
-1 & 0 & 1 & 0 \\
0 & -1 & 0 & 1 
\end{array} \right ) \left ( \begin{array}{c}
P_1 \\
P_2 \\
P_3 \\
P_4
\end{array} \right ) \; \Longleftrightarrow \; \partial  \tau_{(1),i} = \sum_{j} e_{i}^{j} \tau_{(0),j}\;.
\end{equation}
Note that we stretch, twist or bend this little grid without changing the connectivity and orientation, this matrix relation remains valid. This is a metric-free relation because it is independent of shape and size.

Similarly, we can relate the boundary of the surface to the line segments
\begin{equation}
\partial S = \left ( 1\; -1\; -1\; 1\right ) \left ( \begin{array}{c}
L_1 \\
L_2 \\
L_3 \\
L_4
\end{array} \right ) \; \Longleftrightarrow \; \partial  \tau_{(2),i} = \sum_{j} e_{i}^{j} \tau_{(1),j}\;.
\end{equation}
Again, deformations of the grid which do not change the connectivity between surfaces and lines or orientation give the same matrix representation of the boundary operator.

For a chain complex we need to have that $\partial\circ\partial \equiv 0$ and this also follows from the matrix representation of the boundary operator given in this example
\begin{equation}
\partial \partial S = \left ( 1\; -1\; -1\; 1\right ) \left ( \begin{array}{c}
\partial L_1 \\
\partial L_2 \\
\partial L_3 \\
\partial L_4
\end{array} \right ) = \left ( 1\; -1\; -1\; 1\right ) \left ( \begin{array}{cccc}
-1 & 1 & 0 & 0 \\
0 & 0 & -1 & 1 \\
-1 & 0 & 1 & 0 \\
0 & -1 & 0 & 1 
\end{array} \right ) \left ( \begin{array}{c}
P_1 \\
P_2 \\
P_3 \\
P_4
\end{array} \right ) = 0\;.
\end{equation}
\end{example}

Example~\ref{ex:incidence_matrix} shows that on an oriented grid, the boundary operator can be represented by an {\em incidence matrix}, $\incidenceboundary{k}{k-1}$, which relates a $k$-dimensional geometric object to its $(k-1)$-dimensional boundary. We will denote the matrix representing the boundary operator acting on $k$-cells, $\incidenceboundary{k}{k-1}$. The relation $\partial \circ \partial = 0$ then reads $\incidenceboundary{k}{k-1} \incidenceboundary{k-1}{k-2} = 0$.

\begin{remark}\label{rem:incindece_matrix}
The incidence matrices, $\incidenceboundary{k}{k-1}$, are a matrix representation of the boundary operator.
\end{remark}

\begin{remark}
We have that the boundary of the boundary always yields the empty set. If the boundary of a geometric object is the empty set, this does {\em not} necessarily imply that the object itself is a boundary. Only on contractible domains we have that when the boundary of a $k$-dimensional object is empty, that object must be the boundary of a $(k+1)$-dimensional object. This is the Poincar\'{e} lemma.
\end{remark}

Dual to the space of $k$-chains, $C_{k} \left( D \right)$,  is the space of $k$-cochains, $C^{k} \left(D \right)$, defined as the set of all homomorphisms, $\kcochain{c}{k} : C_{k} \left( D \right) \rightarrow \realNumber$
\begin{align}
\left< \kcochain{c}{k}, \kchain{c}{k} \right> := \kcochain{c}{k} \left( \kchain{c}{k} \right).
\end{align}

\begin{remark}
Essentially, a $k$-cochain assigns a physical value to a geometric object. For instance, it can assign mass to a 3-dimensional object, flux to a 2-dimensional object, temperature to a point, or work to a curve. In this respect duality pairing between a cochain and a chain is the discrete analogue of integration in the continuous setting, (\ref{eq:integration_difform})
\end{remark}

\begin{remark}
The value assigned to a $k$-chain (or $k$-cell) is assigned to the geometric object {\em as a whole}. There is no particular point in this object where this value is anchored. So if we assign 1 gram to a 3-dimensional object, 1 gram will be the attribute of that object and we do not assign it to the center of mass. The center of mass is a priori unknown anyway, because we only know the mass of the object and not its mass density distribution.
\end{remark}

\vskip 0.3cm
\noindent
\framebox[1.0\linewidth]{
\begin{minipage}{0.92\linewidth}
With the duality pairing between cochains and chains, one can define the formal adjoint of the boundary operator, the \emph{coboundary operator}, $\delta : C^{k} \left( D \right) \rightarrow C^{k+1} \left( D \right)$,
\begin{align}
\left< \delta \kcochain{c}{k}, \kchain{c}{k+1} \right> := \left< \kcochain{c}{k}, \partial \kchain{c}{k+1} \right>.
\label{Equation::DiscreteStokes}
\end{align}
This equation is the discrete analogue of the generalized Stokes equation, (\ref{eq:generalized_Stokes}).
\end{minipage}}
\vskip 0.3cm

\begin{remark}
The coboundary operator, $\delta$, is the discrete analogue of the exterior derivative, $\ederiv$. Derivatives require a certain smoothness on the objects to be differentiated. Such smoothness assumptions are not available in a discrete, topological setting. By defining the discrete derivative in terms of the boundary operator -- which is well-defined for discrete objects -- we circumvent any smoothness requirements.
\end{remark}

Just like the exterior derivative, the coboundary operator is nilpotent $\delta \delta \kcochain{c}{k} = 0$ for all $\kcochain{c}{k} \in C^{k} \left( D \right)$. This follows directly from the definition and the fact that $\partial \circ \partial \equiv 0$.

Let $C_{k} \left( D \right)$ be the space of $k$-chains with basis $\left\{ \tau_{(k),j} \right\}$, then a dual basis, $\left\{ \tau^{(k),i} \right\}$ of $C^{k} \left( D \right)$ is given, such that $\tau^{(k),i}  \left( \tau_{(k),j}  \right) = \delta_{j}^{i}$. All $k$-cochains can be represented as linear combinations of these basis elements,
\begin{align}
\kcochain{c}{k} = \sum_{i}c_{i} \tau^{(k),i}.
\label{eq:expansion_k-cochain}
\end{align}

\begin{remark}
Note also here the resemblance between $\tau^{(k),i}  \left( \tau_{(k),j}  \right) = \delta_{j}^{i}$ and $\ederiv x^i ( \partial/\partial x^j) = \delta_j^i$ for the duality pairing between vectors and forms. The $k$-cochains therefore acts as the discrete differential $k$-form.
\end{remark}

Once a basis for $k$-cochains has been chosen, any $k$-cochain is completely determined by the expansion coefficients, $c_i$, in (\ref{eq:expansion_k-cochain}). That is, there is an isomorphism $\psi$
\[ \psi\,:\, C^k(D) \mapsto \mathbb{R}^p\;,\;\; p = \mbox{rank}(C^k(D)) = \mbox{rank}(C_k(D)) \;,\]
defined by \label{page:isomorphism}
\begin{equation}
\psi(\kcochain{c}{k}) = \psi \left ( \sum_{i}c_{i} \tau^{(k),i} \right ) = [ c_1\,\dots\,c_p ], \;\; p = \mbox{rank}(C^k(D)) \;.
\label{eq:isomorphism_cochains_expansion_coeff}
\end{equation}
Note that a $k$-cochain is represented by a {\em row vector}, instead of a column vector for the expansion coefficients for the chains. Duality pairing of a cochain and chain in terms of the expansion coefficients then simply reduces to
\[ \left< \kcochain{c}{k}, \kchain{c}{k} \right> = \sum_{i=1}^p c_i c^i \;.\]
Again note the resemblance with the duality paring at the continuous level given by (\ref{eq:duality_pairing_vector_covector}).
 
With the isomorphism which identifies the chains and cochains with their expansion coefficients, we also have a natural matrix representation for the coboundary operator \label{page:coboundary_incidence}

\begin{eqnarray*}
\left< \kcochain{c}{k}, \partial \kchain{c}{k+1} \right> & = & \sum_{i=1}^{\mbox{rank}(C_k(D))} c_i \left ( \incidenceboundary{k+1}{k} c^i \right ) \\
 & = & \sum_{i=1}^{\mbox{rank}(C_k(D))} \left ( c_i \incidenceboundary{k+1}{k} \right ) c^i  \\
 & = & \left< \delta \kcochain{c}{k}, \kchain{c}{k+1} \right> \;.
\end{eqnarray*}
So, when the row vector $c^{(k)}$ contains the expansion coefficients, $c_i$, for the cochain $\kcochain{c}{k}$, then the row vector $c^{(k)}\incidenceboundary{k+1}{k}$ contains the expansion coefficients for $\delta \kcochain{c}{k}$.

\begin{remark}
This relation is one of the remarkable properties of mimetic methods. Bear in mind that the coboundary operator encodes the action of the gradient ($k=0$), the curl ($k=1$) and the divergence operator ($k=2$) at the discrete level. Once a basis has been chosen, the matrix operation which represents this operation is given by the incidence matrix of the oriented grid. These relations are {\em exact}, no approximations have been performed. So the topology of the oriented grid determines the action of the grad, curl and div.
\end{remark}

In Section~\ref{sec:Orientation} we made a distinction between true forms and pseudo-forms. The operator which switches between the two is the Hodge-$\star$ operator which was defined in (\ref{eq:SpacesHodge}) and (\ref{eq:HodgeInnerProduct}). This cannot be accomplished on a single cell complex. Consider Example~\ref{ex:incidence_matrix}. This is a 2-dimensional problem, therefore $n=2$. The Hodge-$\star$ is a bijection between $k$-forms and $(n-k)$-forms. In the discrete setting of Example~\ref{ex:incidence_matrix}, it needs to be a bijection between $k$-cochains and $(n-k)$-cochains. However, in this example we have four $0$-cells (points) and only one $(n-k)=(2-0)$-cells (surfaces). Besides, not only the number of dual cells differs from the number of associated primal cells, but also the kind of integrals we want to represent differs, see Section~\ref{sec:Orientation}. So there cannot be a bijection between $0$-cochains and $2$-cochains. This is not a particular problem for the grid presented in Example~\ref{ex:incidence_matrix}, but for grids (cell complexes), in general.

\begin{figure}[ht!]
\centering
\includegraphics[width=0.75\textwidth]{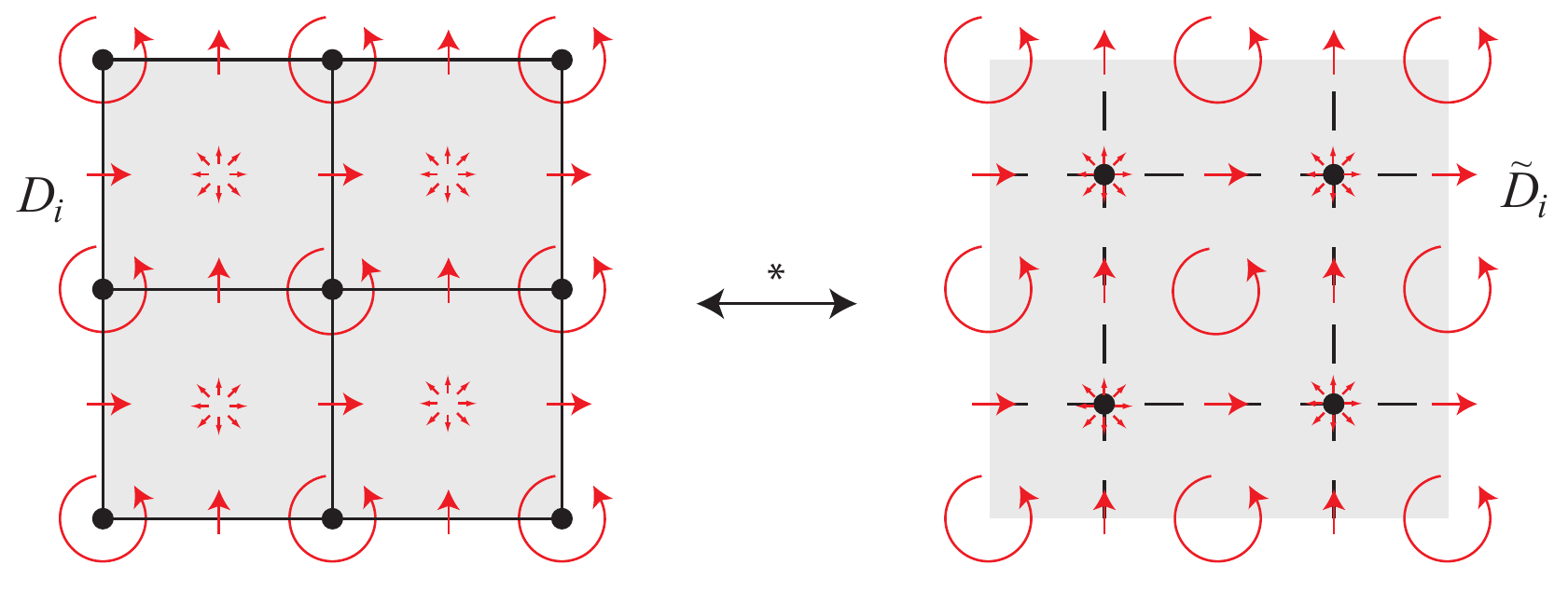}
\caption{A cell complex (left) and its associated dual grid (right)}
\label{fig:dualcellcomplex2D_top}
\end{figure}

In order to incorporate the duality between of $k$- and $(n-k)$-forms at the discrete level we construct a {\em dual grid}\label{page:dual_grid}. The construction is as follows: With every $k$-cell in the cell complex, $D$, we associate an $(n-k)$-cell in the dual complex. An example of such a construction is shown in Figure~\ref{fig:dualcellcomplex2D_top}. The collection of all dual cells constitutes the dual grid, $\tilde{D}_i$, which contains the interior dual cells of the dual cell complex $\tilde{D}$\label{page:tilde_D}, see \citep{kreeft2011mimetic} for further details. 
In Figure~\ref{fig:dualcellcomplex2D_top} the pseudo-forms (outer-oriented forms) are represented on the cell complex $D$ on the left and the true forms (inner-oriented forms) are represented on the dual grid. Note that the dual grid itself is not a cell complex in general, because for a cell complex we require that if an $k$-dimensional element is part of the grid, then so is its boundary. That is not the case for the dual grid. There are some line segments for which not all boundary points are part of the grid. We will see in the remainder of this paper that this 'missing boundary' for the dual grid leads to so-called ghost points in finite volume methods and to boundary integrals in finite element methods.

%

\section{SWITCHING BETWEEN CONTINUOUS AND DISCRETE}
\label{Section::Discretization}
In the previous section we described integrals and integral relations at the continuous level in terms of differential geometry and at the discrete level in terms of chains and cochains. In this section we explain how to convert the continuous forms into discrete cochains and vice versa.

In this paper we focus on the Laplace operator given by \eqref{eq:Laplace-DeRham}. In Section~\ref{Section::DifferentialGeometry} we presented this operator at the continuous level in terms of differential forms. In Section~\ref{Section::AlgebraicTopology} a discrete description was given. In this section we are going to introduce the operations which will allow us to switch between the continuous formulation and the discrete formulation.

\subsection{From continuous to discrete}

The reduction operator, $\reduction : \differentialFormSpace{k}{\manifold{M}} \rightarrow C^{k} \left( D \right)$ maps differential forms to cochains by 
\begin{align}
\left< \reduction \alpha^{(k)}, \tau_{(k),i} \right> := \int_{\tau_{(k),i}} \alpha^{(k)} = \left< \kdifform{\alpha}{k},\tau_{(k),i} \right>.
\end{align}

Then, for all $\kchain{c}{k} \in C_{k}\left(D\right)$, the reduction of the $k$-form, $\alpha^{(k)} \in \differentialFormSpace{k}{\manifold{M}}$, to the $k$-cochain, $a^{(k)} \in C^{k} \left(D \right)$ is,

\begin{align}
a^{(k)} \left( \kchain{c}{k} \right) := \left< \reduction \alpha^{(k)}, \kchain{c}{k} \right> = \sum_{i} c^{i} \left< \reduction \alpha^{(k)}, \tau_{(k),i} \right> = \sum_{i} c^{i} \int_{\tau_{(k),i}} \reduction \alpha^{(k)} = \int_{\kchain{c}{k}} \alpha^{(k)}.
\end{align}

The reduction has an important commuting property with respect to differentiation in terms of exterior derivative and coboundary operator,
\begin{align}
\reduction d = \delta \reduction.
\end{align}

\begin{remark}
Commuting relations are essential to mimetic methods. Essentially they state that operations at the continuous level are mimicked by equivalent relations at the discrete level. In this case, it makes no difference whether we take the derivative and then convert to discrete variables or first map to discrete variables and then take the discrete derivative.
\end{remark}

\subsection{From discrete to continuous}

The other crucial operator is the reconstruction map, $\reconstruction : C^{k} \left( D \right) \longmapsto \differentialFormSpace{k}{\manifold{M}} $. This operator needs to have the following commuting property,
\begin{align}
d \reconstruction = \reconstruction \delta \;.
\label{eq:commuting_reconstruction_derivative}
\end{align}
The map $\reconstruction$ is an injective, but not surjective.

\begin{remark}
The reconstruction operator is essentially an interpolation operator which interpolates the discrete cochains to continuous $k$-forms. The commuting property (\ref{eq:commuting_reconstruction_derivative}) states that we can reconstruct and take the derivative or first take the discrete derivative and then reconstruct. It should give the same differential form. This condition poses severe restriction on the admissible reconstruction functions. For triangular elements these reconstructions are known as Whitney forms, \citep{bossavit:japanese_01}.
\end{remark}

Furthermore, the reconstruction operator must be the right inverse of $\reduction$,  so $\reduction \reconstruction = Id$ on $C^{k} \left( D \right)$ and it should approximate the left inverse of $\reduction$, so $\reconstruction \reduction = Id + {\mathcal O} \left( h^{p} \right)$.

\begin{remark}
The first condition, $\reduction \reconstruction = Id$, is a {\em consistency condition}. The second condition, $\reconstruction \reduction = Id + {\mathcal O} \left( h^{p} \right)$, is an {\em approximability condition}. Note that the $O \left( h^{p} \right)$-term is in the kernel of the reduction operator, i.e. $\reduction {\mathcal O} \left( h^{p} \right) = 0$.
\end{remark} 

\subsection{Mimetic projection}

The mimetic projection, $\projection$,  is given by 
\begin{equation} 
\projection := \reconstruction \circ \reduction \;.
\label{eq:mimetic_projection}
\end{equation}
\begin{equation*}
\begin{tikzpicture}[descr/.style={fill=white,inner sep=2.5pt}]
\matrix (m) [matrix of math nodes, row sep=1.5em,
column sep=1.5em]
{\differentialFormSpace{k}{\manifold{M}} & \differentialFormSpaceDiscrete{k}{\manifold{M}; C_{k}} \\
C^{k} \left( D \right) &   \\};
\path[->,font=\scriptsize]
(m-1-1) edge node[auto] {$\projection$} (m-1-2)
(m-1-1) edge node[left] {$\reduction$} (m-2-1)
(m-2-1) edge node[below] {$\reconstruction$} (m-1-2);
\end{tikzpicture}
\end{equation*}
where $\differentialFormSpaceDiscrete{k}{\manifold{M}; C_{k}}$ denotes the range of of the reconstruction, $\reconstruction$, of $k$-cochains associated with space of $k$-chains,


\vskip 0.3cm
\noindent
\framebox[1.0\linewidth]{
\begin{minipage}{0.92\linewidth}
Due to the commuting properties of reduction and reconstruction we have
\begin{equation} 
\projection \ederiv = \reconstruction \reduction \ederiv = \reconstruction \delta \reduction  = \ederiv \reconstruction \reduction  = \ederiv \projection \;.
\label{eq:commutation_relations_projection_derivative}
\end{equation}
\end{minipage}}
\vskip 0.3cm

\begin{remark}
This commuting property is the most important one of all and it has important consequences. For instance, let $\kdifform{\alpha}{k}$ be in the null space of $\ederiv$, i.e. $\ederiv \kdifform{\alpha}{k} = \kdifform{0}{k+1}$, then its projected form, $\kdifformh{\alpha}{k}$, will also be in the null space of $\ederiv$. Conforming null spaces play an important role in mixed finite element methods, \citep{Kreeft_A_priori_Stokes,brezzi1991mixed}, and the existence of discrete potentials, \citep{RobidouxSteinberg2011}.
\end{remark}

\vskip 0.3cm
\noindent
\framebox[1.0\linewidth]{
\begin{minipage}{0.92\linewidth}
Another important property is the commuting relation between the mimetic projection, $\projection$ and the pullback of a map $\Phi^\star$
\begin{equation} 
\projection \Phi^\star = \Phi^\star \projection \;.
\label{eq:commuting_projection_pullback}
\end{equation}
This will allow us to do all computations on a reference element instead of working directly in curvilinear coordinates.
For a proof of this property, see \citep{kreeft2011mimetic}.
\end{minipage}}
\vskip 0.3cm

Having defined the reduction and reconstruction operators for the cell complex a similar derivation can be done for the dual cell complex, $\tilde{D}$, see its construction on page~\pageref{page:tilde_D}.

The space $\differentialFormSpaceDiscrete{k}{\manifold{M}; \tilde{C}_{k}}$ is the space of discrete $k$-forms, with $k$-cochains associated with the dual $k$-cells.

\begin{align}
 \differentialFormSpaceDiscrete{k}{\manifold{M}; \tilde{C}_{k}} = \dualprojection\differentialFormSpaceDual{k}{\manifold{M}} := \reconstructionDual \reductionDual \differentialFormSpace{k}{\manifold{M}} \;,
\end{align}
where $\tilde{C}_k$ is the space of $k$-chains on the dual complex $\tilde{D}$, $\reductionDual$, is the reduction of differential forms on the dual chains and $\reconstructionDual$ constitutes the reconstruction of differential forms from cochains defined on the dual complex.


A discrete wedge product is introduced such that $\wedge_h:\Lambda^k_h\times\Lambda^l_h\rightarrow\Lambda^{k+l}_h$, given by
\begin{equation}
\kdifformh{\a}{k}\wedge_h \kdifformh{\b}{l}:=\pi\left(\kdifformh{\a}{k}\wedge\kdifformh{\b}{l}\right)\;.
\label{reductionwedge}
\end{equation}

\begin{remark}
The discrete wedge product satisfies all properties of the continuous wedge product except associativity. Associativity can be restored by using a projection $\pi_{\bar{h}}$, $\bar{h}<h$, i.e. a mimetic projection on a refined cell complex. The precise refinement needed to fully represent the wedge product at the finite dimensional level depends on the number of terms in the wedge product and the type of reconstruction forms.
\end{remark}

%

There are essentially two different ways in which we can incorporate the action of the Hodge-$\star$ operator in the finite dimensional setting. Either we can use the definition of the inner product for differential forms, (\ref{eq:HodgeInnerProduct}) or we can use reconstruction of the cochains to obtain differential form, then apply the Hodge-$\star$ operator and subsequently we reduce the result on the topological dual grid, i.e. $\tilde{\reduction}\star\reconstruction$ or $\reduction\star\tilde{\reconstruction}$.
In the former approach we make use of
\begin{equation}
\innerspace{\kdifformh{\a}{k}}{\kdifformh{\b}{k}}{L^2\manifold{M}}=\int_\manifold{M} \kdifformh{\a}{k}\wedge\star \kdifformh{\b}{k}.
\end{equation}
In the definition of this inner product, the Hodge is taken of the second argument and therefore the inner product implicitly. Methods based on the use of this inner product will be referred to as the {\em single grid method}.

The method where the Hodge is applied to reconstructed cochains and reduced onto the dual grid will be called {\em dual grid method}, because for this approach an explicit dual grid needs to be defined.

%

\subsection{Basis forms}

Let $\manifold{M}$ be the computational domain decomposed into $M$ non-overlapping, possibly curvilinear quadrilateral or hexahedral closed sub-domains, $Q_{m}$,

\begin{align}
\manifold{M}  = \bigcup_{m=1}^{M} Q_{m}, \quad Q_{m} \cap Q_{l} = \partial Q_{m} \cap \partial Q_{l}, \quad m \neq l,
\end{align}

where in each sub-domain a Gauss-Lobatto mesh is constructed. The collection of Gauss-Lobatto grids in all elements constitutes the cell complex $D$. For each element $Q_{m}$ there exists a sub cell complex, $D_{m}$. Note that $D_{m} \cap D_{l}$, $m \neq l$, is not an empty set in the case they are neighboring elements, but contains all $k$-cells, $k < n $, of the common boundary. 

Each sub-domain is a map from the reference element $Q_{ref} = \left[-1, 1\right]^{n}$, $n = \mbox{dim} \left( \manifold{M} \right)$ using the mapping $\Phi_{m}: Q_{ref} \rightarrow Q_{m}$. All differential forms defined on $\manifold{M}_{m}$ are pulled back onto the reference element using the following pullback operation $\pullback_{m}: \differentialFormSpaceDiscrete{k}{Q_{m}, C_{k}} \rightarrow \differentialFormSpaceDiscrete{k}{Q_{ref}, C_{k}}$.

The cochains are approximated using piecewise polynomial expansions on the quadrilateral or hexahedral elements using tensor products. Thus, it suffices to derive the basis forms in one dimension and afterwards construct the $n$-dimensional basis forms. Furthermore, because of the commutation between the projection operator and the pullback, (\ref{eq:commuting_projection_pullback}), only the interpolation for the reference element is shown. 

Consider a 0-form, $\alpha^{(0)} \in \differentialFormSpace{0}{Q_{ref}}$, where $Q_{ref} := \xi \in \left[ -1,1 \right]$, on which a cell complex $D$ consists on $N+1$ nodes $\xi_{i}$, where $-1 \leq \xi_{0} < \cdots < \xi_{N} \leq 1$, and $N$ edges, $\tau_{(1),i} = \left[ \xi_{i-1},\xi_{i} \right]$, of which the nodes are the boundaries. Corresponding to this set of nodes (0-chain) there exists a projection,$\pi_{h}$, using the $N^{th}$ order \emph{Lagrange polynomials}, $h_{i} \left( \xi \right)$, to approximate a 0-form,
\begin{align}
\pi_{h} \alpha^{(0)} \left( \xi \right) = \sum_{i=0}^{N} a_{i} h_{i} \left( \xi \right).
\end{align}

\begin{figure}[ht!]
\centering
\includegraphics[width=0.40\textwidth]{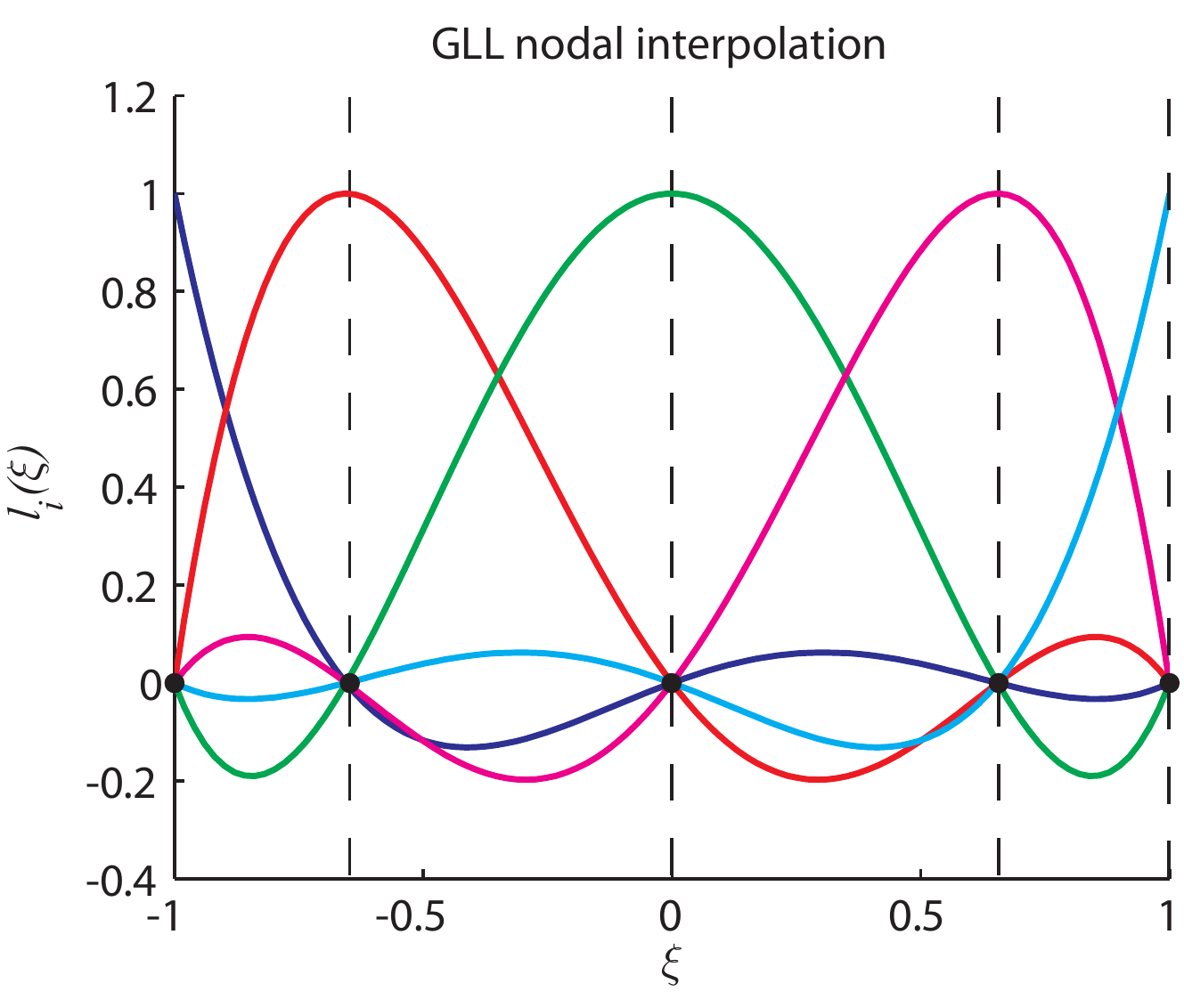}
\caption{The Lagrange polynomials, $h_i(\xi)$, for $N=4$}
\label{fig:gll}
\end{figure}

Lagrange polynomials interpolate nodal values. Thus, they are suitable to reconstruct the cochain $a^{(0)} = \reduction \alpha^{(0)}$. These polynomials are constructed such that their value is one in the corresponding point and zero in all other mesh points, see Figure~\ref{fig:gll}
\begin{align}
\reduction h_{i}^{(0)} \left( \xi \right) = h_{i}^{(0)} \left(\xi_{p} \right) = \begin{cases}
1 & \text{if } i = p \\
0 & \text{if } i \neq p \\
\end{cases}.
\end{align} 

Similarly for the projection of 1-forms Gerritsma \citep{gerritsma::edge_basis} and Robidoux \citep{robidoux-polynomial} derived 1-form polynomials called \emph{edge polynomials}, $e_{i}\left( \xi \right) \in \differentialFormSpaceDiscrete{1}{Q_{ref};C_{1}}$. More details can be found in \citep{kreeft2011mimetic}.
\begin{align}
e_{i}\left( \xi \right) = \epsilon_{i} \left( \xi \right) d\xi, \quad \text{with} \quad \epsilon_{i} = - \sum_{k=0}^{i-1} \frac{d h_{k}}{d \xi}.
\label{eq:edge_functions_e_and_epsilon}
\end{align}
The cochain corresponding to the line segment (1-cell), $\tau_{(1),i}$ is given by $u_{i} = a_{i} - a_{i-1}$ and so $u^{(1)} = \delta a^{(0)}$ is the discrete derivative operator in 1D. This is a purely topological operation and $d \reconstruction a^{(0)} = \reconstruction \delta a^{(0)}$. The 1-form edge polynomial can be separated into its polynomial and its basis,

\begin{figure}[ht!]
\centering
\includegraphics[width=0.40\textwidth]{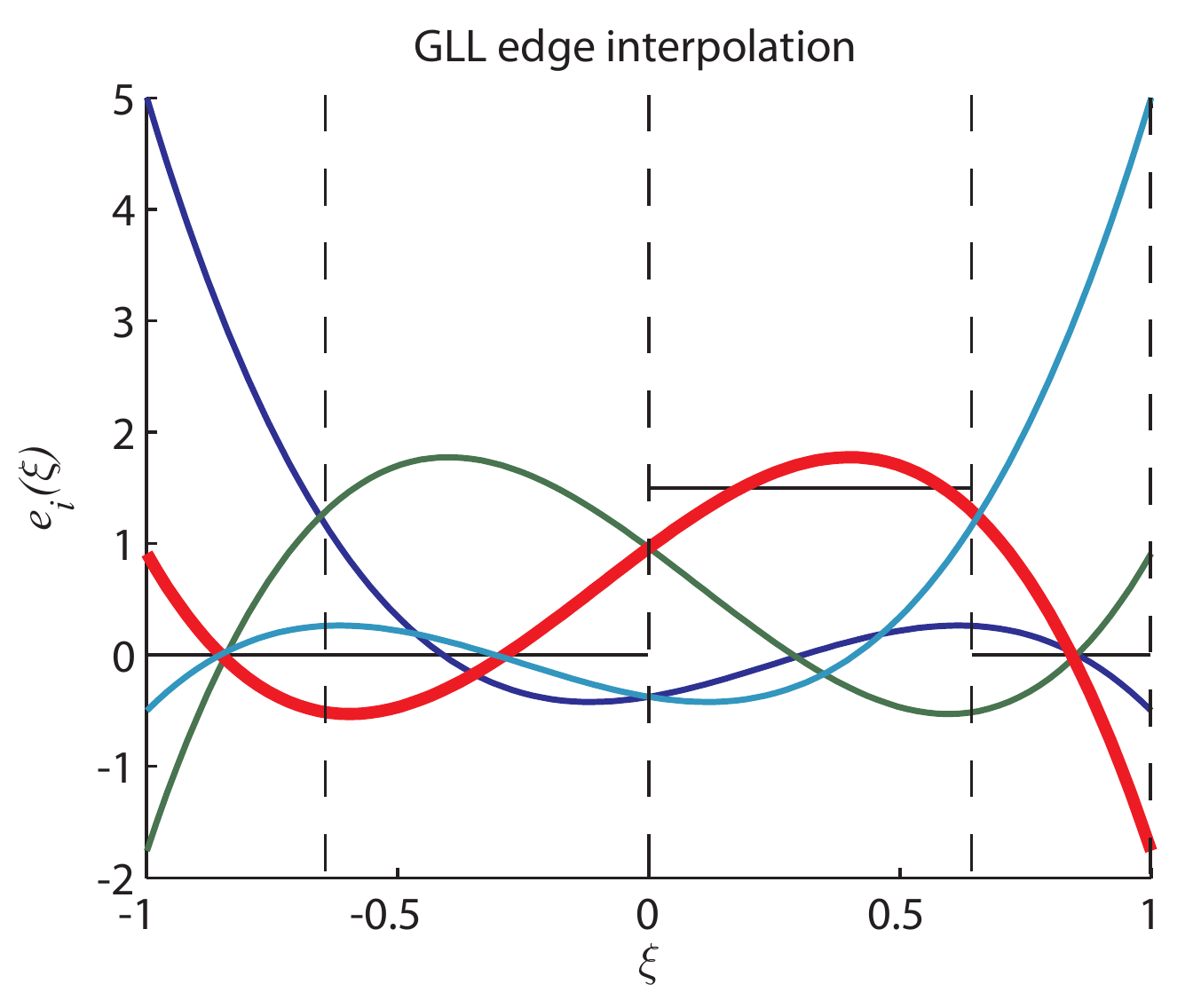}
\caption{Edge polynomial $1$-forms, $e_i(\xi)$, for $N=4$}
\label{fig:edge}
\end{figure}

The edge functions are constructed such that when integrating $e_{i} \left( \xi \right)$ over a line segment it gives one for the corresponding element and zero for any other line segment, see Figure~\ref{fig:edge}
\begin{align}
\reduction e_{i} \left( \xi \right)  = \int_{\xi_{p-1}}^{\xi_{p}} e_{i} \left( \xi \right) = \begin{cases}
1 & \text{if } i = p \\
0 & \text{if } i \neq p \\
\end{cases}
\end{align}

The full connection between continuous and discrete counterparts can be summarized in the following diagram,
\begin{center}
\begin{tikzpicture}
  \matrix (m) [matrix of math nodes, row sep=1em,
    column sep=1em]{
    & \cochainSpace{0}{D} & & \cochainSpace{1}{D} & & \cochainSpace{2}{D} & & \cochainSpace{3}{D}  \\
    \differentialFormSpace{0}{\manifold{M}} & & \differentialFormSpace{1}{\manifold{M}}  & & \differentialFormSpace{2}{\manifold{M}}  & & \differentialFormSpace{3}{\manifold{M}} &\\
    & \cochainSpaceDual{3}{D} & & \cochainSpaceDual{2}{D} & & \cochainSpaceDual{1}{D} & & \cochainSpaceDual{0}{D} \\
    \differentialFormSpaceDual{3}{\manifold{M}} & & \differentialFormSpaceDual{2}{\manifold{M}}& & \differentialFormSpaceDual{1}{\manifold{M}} & & \differentialFormSpaceDual{0}{\manifold{M}} \\};
  \path[-stealth]
    (m-1-2) edge node[below,xshift=-0.3cm] {$\star_{h}$} (m-3-2)
    (m-1-4) edge node[below,xshift=0.3cm] {$\star_{h}$} (m-3-4)
    (m-1-6) edge node[below,xshift=0.3cm] {$\star_{h}$} (m-3-6)
    (m-1-8) edge node[below,xshift=0.3cm] {$\star_{h}$} (m-3-8)

    (m-3-2) edge node[below,xshift=-0.3cm] {} (m-1-2)
    (m-3-4) edge node[below,xshift=0.3cm] {} (m-1-4)
    (m-3-6) edge node[below,xshift=0.3cm] {} (m-1-6)
    (m-3-8) edge node[below,xshift=0.3cm] {} (m-1-8)
    (m-2-1) edge node[below,xshift=-0.3cm,yshift=0.5cm] {$\star$} (m-4-1)
    (m-2-3) edge node[below,xshift=0.3cm,yshift=0.5cm] {$\star$} (m-4-3)
    (m-2-5) edge node[below,xshift=0.3cm,yshift=0.5cm] {$\star$} (m-4-5)
    (m-2-7) edge node[below,xshift=0.3cm,yshift=0.5cm] {$\star$} (m-4-7)

    (m-4-1) edge node[below,xshift=-0.3cm,yshift=0.5cm] {} (m-2-1)
    (m-4-3) edge node[below,xshift=0.3cm,yshift=0.5cm] {} (m-2-3)
    (m-4-5) edge node[below,xshift=0.3cm,yshift=0.5cm] {} (m-2-5)
    (m-4-7) edge node[below,xshift=0.3cm,yshift=0.5cm] {} (m-2-7)
    (m-2-1) edge node[above,xshift=0.3cm] {$d$} (m-2-3)
    (m-2-3) edge node[above,xshift=0.3cm] {$d$} (m-2-5)
    (m-2-5) edge node[above,xshift=0.3cm] {$d$} (m-2-7)

    (m-4-3) edge node[below] {$d$} (m-4-1)
    (m-4-5) edge node[below] {$d$} (m-4-3)
    (m-4-7) edge node[below] {$d$} (m-4-5)
    (m-1-2) edge node[above,xshift=-0.3cm] {$\delta$} (m-1-4)
    (m-1-4) edge node[above,xshift=-0.3cm] {$\delta$} (m-1-6)
    (m-1-6) edge node[above,xshift=-0.3cm] {$\delta$} (m-1-8)

    (m-3-4) edge node[below,xshift=-0.3cm] {$\delta$} (m-3-2)
    (m-3-6) edge node[below,xshift=-0.3cm] {$\delta$} (m-3-4)
    (m-3-8) edge node[below,xshift=-0.3cm] {$\delta$} (m-3-6)
    (m-2-1) edge[bend angle=15, bend right] node[right] {$\reduction^0$} (m-1-2)
    (m-2-3) edge[bend angle=15, bend right] node[right] {$\reduction^1$} (m-1-4)
    (m-2-5) edge[bend angle=15, bend right] node[right] {$\reduction^2$} (m-1-6)
    (m-2-7) edge[bend angle=15, bend right] node[right] {$\reduction^3$} (m-1-8)

    (m-4-1) edge[bend angle=15, bend right] node[right,xshift=0.1cm] {$\tilde{\reduction}^3$} (m-3-2)
    (m-4-3) edge[bend angle=15, bend right] node[right,xshift=0.1cm] {$\tilde{\reduction}^2$} (m-3-4)
    (m-4-5) edge[bend angle=15, bend right] node[right,xshift=0.1cm] {$\tilde{\reduction}^1$} (m-3-6)
    (m-4-7) edge[bend angle=15, bend right] node[right,xshift=0.1cm] {$\tilde{\reduction}^0$} (m-3-8)
    
    (m-1-2) edge[bend angle=15, bend right] node[left] {$\reconstruction^0$} (m-2-1)
    (m-1-4) edge[bend angle=15, bend right] node[left] {$\reconstruction^1$} (m-2-3)
    (m-1-6) edge[bend angle=15, bend right] node[left] {$\reconstruction^2$} (m-2-5)
    (m-1-8) edge[bend angle=15, bend right] node[left] {$\reconstruction^3$} (m-2-7)

    (m-3-2) edge[bend angle=15, bend right] node[left] {$\tilde{\reconstruction}^3$} (m-4-1)
    (m-3-4) edge[bend angle=15, bend right] node[left] {$\tilde{\reconstruction}^2$} (m-4-3)
    (m-3-6) edge[bend angle=15, bend right] node[left] {$\tilde{\reconstruction}^1$} (m-4-5)
    (m-3-8) edge[bend angle=15, bend right] node[left] {$\tilde{\reconstruction}^0$} (m-4-7);
\end{tikzpicture}
\end{center}


The mimetic framework uses Lagrange, $h_i(\xi)\in H\Lambda^0(\manifold{M}_{\rm ref})$, and edge functions, $e_i(\xi)\in L^2\Lambda^1(\manifold{M}_{\rm ref})$, for the reconstruction, $\mathcal{I}$, where the latter is constructed using the former; i.e., from the finite dimensional 0-form $\pi_ha^{(0)}=\sum_{i=0}^N a_ih_i(\xi)\in\Lambda^0_h(\manifold{M}_{\rm ref};C_0)$, we define $\pi_hb^{(1)}\in\Lambda^1_h(\manifold{M}_{\rm ref};C_1)$, such that
\[
\pi_hb^{(1)}=\pi_h\ederiv a^{(0)}=\sum_{i=1}^N b_i e_i(\xi)\;,
\]
where $b_i=a_i-a_{i-1}$. Because we consider tensor products to construct higher-dimensional interpolation, it is sufficient to show that the projection operator is bounded in one dimension, \citep{kreeft2011mimetic}.

All the interpolation functions were defined in one dimension the extension to multidimensional is straightforward by means of tensor products,
\begin{align}
 \begin{array}{l}
 \kdifform{P_{i,j,k}}{0} \left( \xi, \eta, \zeta \right) = h_{i}\left(\xi \right) \otimes h_{j}\left(\eta \right) \otimes h_{k}\left(\zeta \right) \\
 \kdifform{L_{i,j,k}}{1} \left( \xi, \eta, \zeta \right) = \left\{ e_{i}\left(\xi \right) \otimes h_{j}\left(\eta \right) \otimes h_{k}\left(\zeta \right), h_{i}\left(\xi \right) \otimes e_{j}\left(\eta \right) \otimes h_{k}\left(\zeta \right), h_{i}\left(\xi \right) \otimes h_{j}\left(\eta \right) \otimes e_{k}\left(\zeta \right)  \right\} \\
 \kdifform{S_{i,j,k}}{2} \left( \xi, \eta, \zeta \right) = \left\{ h_{i}\left(\xi \right) \otimes e_{j}\left(\eta \right) \otimes e_{k}\left(\zeta \right), e_{i}\left(\xi \right) \otimes h_{j}\left(\eta \right) \otimes e_{k}\left(\zeta \right), e_{i}\left(\xi \right) \otimes e_{j}\left(\eta \right) \otimes h_{k}\left(\zeta \right)  \right\} \\
 \kdifform{V_{i,j,k}}{3} \left( \xi, \eta, \zeta \right) = e_{i}\left(\xi \right) \otimes e_{j}\left(\eta \right) \otimes e_{k}\left(\zeta \right). 
 \end{array} 
 \label{eq:basis_k_forms}
\end{align}

The approximation spaces are spanned by combinations of Lagrange and edge basis forms given by,
\begin{align}
\begin{array}{l}
 \Lambda^{0}_{h}\left(\manifold{M}_{m}\right) = \text{span}\left\{ \kdifform{P_{i,j,k}}{0} \right\}_{i=0, j=0, k=0}^{N,N,N} \\
 \Lambda^{1}_{h}\left(\manifold{M}_{m}\right) = \text{span}\left\{ \left(\kdifform{L_{i,j,k}}{1}\right)_{1} \right\}_{i=1, j=0, k=0}^{N,N,N} \times \text{span}\left\{ \left(\kdifform{L_{i,j,k}}{1}\right)_{2} \right\}_{i=0, j=1, k=0}^{N,N,N} \times \text{span}\left\{ \left(\kdifform{L_{i,j,k}}{1}\right)_{3} \right\}_{i=0, j=0, k=1}^{N,N,N} \\
 \Lambda^{2}_{h}\left(\manifold{M}_{m}\right) = \text{span}\left\{ \left(\kdifform{S_{i,j,k}}{2}\right)_{1} \right\}_{i=0, j=1, k=1}^{N,N,N} \times \text{span}\left\{ \left(\kdifform{S_{i,j,k}}{2}\right)_{2} \right\}_{i=1, j=0, k=1}^{N,N,N} \times \text{span}\left\{ \left(\kdifform{S_{i,j,k}}{2}\right)_{3} \right\}_{i=1, j=1, k=0}^{N,N,N} \\
 \Lambda^{3}_{h}\left(\manifold{M}_{m}\right) = \text{span}\left\{ \left(\kdifform{V_{i,j,k}}{3}\right)_{1} \right\}_{i=1, j=1, k=1}^{N,N,N}
\end{array}
\end{align}

Using the Lagrange polynomials and the associated edge polynomials on the dual grid and applying a tensor product construction, basis forms for the dual complex is constructed similarly.

\section{POISSON EQUATION FOR VOLUME FORMS}
\label{sec:Poisson_volume_forms}

In this section we want to illustrate how both the dual grid and the single grid approach can be used in practice. As a test problem we take the Poisson equation for a $n$-form. The dual grid approach will resemble a staggered finite volume method, whereas the single grid approach leads to a stable, well-posed mixed finite element formulation. 

The Laplace operator was defined in (\ref{eq:Laplace-DeRham}). Since $\ederiv \kdifform{\omega}{n}\equiv 0$ for an $n$-form, the Poisson equation for an $n$-form reads
\[ \Delta \kdifform{\omega}{n} = \kdifform{f}{n}  \quad \Longleftrightarrow \quad \ederiv \ederiv^* \kdifform{\omega}{n} = \kdifform{f}{n}\;.\]
Here we present results for $n=2$ but the method can be readily extended to other values of $n$.

\subsection{Dual grid approach}
\begin{figure}[ht!]
	\centering
		\includegraphics[width=0.75\textwidth]{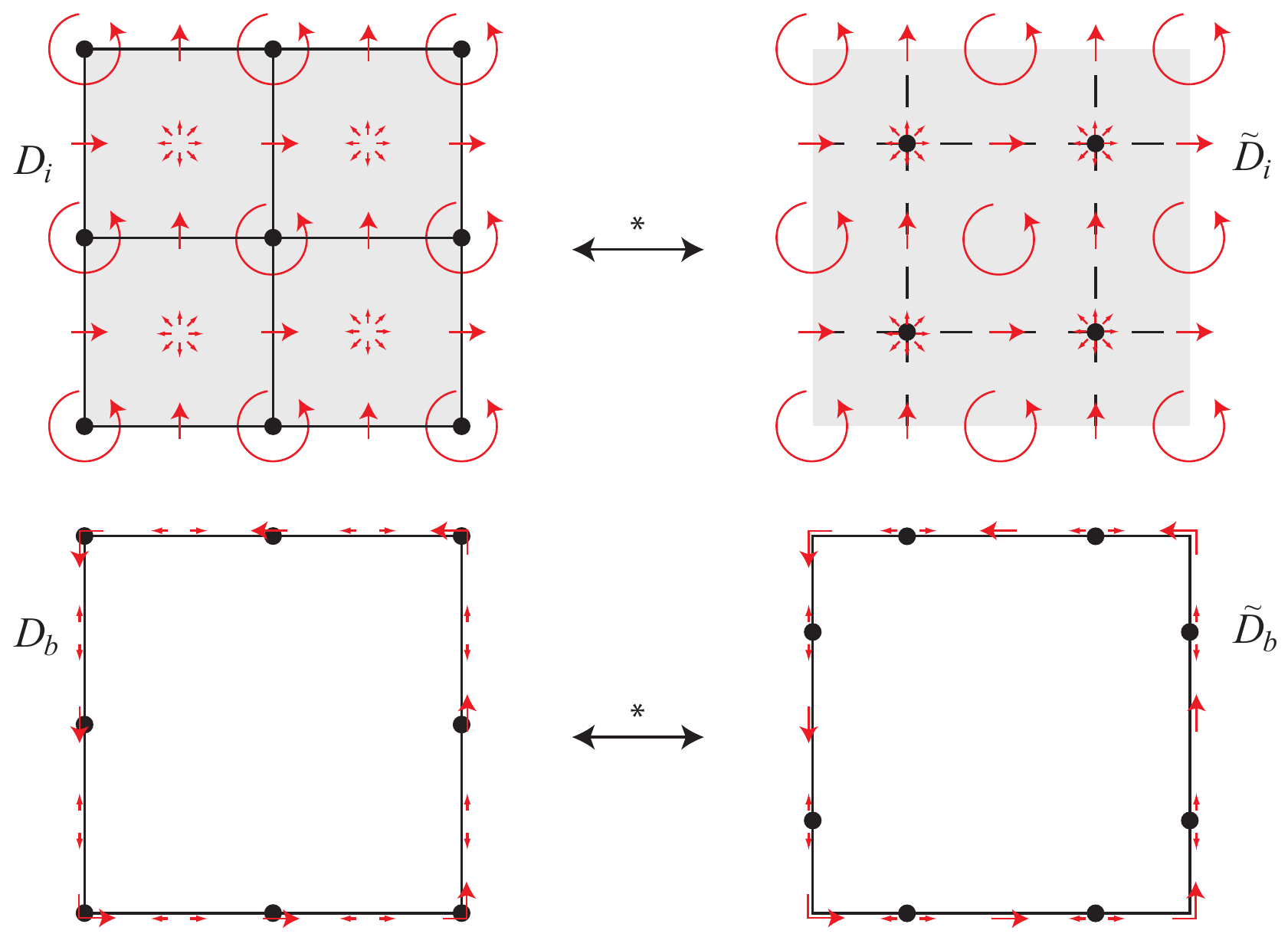}
	\caption{Primal cell complex (top left) and its associated dual grid (top right). Since the dual grid is not a cell complex, we complement this grid with the minimal number of $k$-chains to complete the complex. This completion is the boundary of the dual complex (bottom right). The dual boundary is itself a complex and its dual is the primal boundary complex (bottom left)}
	\label{fig:mesh}
\end{figure}

In the dual grid approach we make use of the fact that $\ederiv^* = (-1)^{n(k+1)+1}\star \ederiv \star$, which for $n=2$ always yields $\ederiv^*=-\star \ederiv \star$. In the dual grid approach we make use of two dual grids as discussed at the end of Section~\ref{Section::AlgebraicTopology} on page~\pageref{page:dual_grid}. We assume that any quadrilateral element, $\manifold{M}_m$, in the $(x,y)$-plane is obtained from a map $\Phi\,:\,(\xi,\eta)\in [-1,1]^2 \rightarrow (x,y)\in \manifold{M}_m$. Then the pullback $\Phi^\star$ maps forms in physical space, $\manifold{M}_m$, to forms on the reference element $[-1,1]^2$. It therefore suffices to explore the analysis on the reference domain. The $2$-forms $\kdifform{\omega}{n}$ and $\kdifform{f}{n}$ will both be represented on the primal grid, i.e. on the cell complex shown top left in Figure~\ref{fig:mesh}. In terms of the basis functions introduced in Section~\ref{Section::Discretization}, these forms are expanded as
\[ \kdifformh{\omega}{2}=\sum_{i=1}^N \sum_{j=1}^N \omega_{ij}e_i(\xi)e_j(\eta) \quad \mbox{ and } \quad \kdifformh{f}{2}=\sum_{i=1}^N \sum_{j=1}^N f_{ij}e_i(\xi)e_j(\eta) \;,\]
where $N$ denotes the number of surfaces in the $x$- and $y$-direction in primal grid of Figure~\ref{fig:mesh}. Here the coefficients $\omega_{ij}$ are the unknown coefficients and $f_{ij}$ is given by
\[ f_{ij} = \int_{\xi_{i-1}}^{\xi_i}\int_{\eta_{j-1}}^{\eta_j} \Phi^\star \kdifform{f}{2} \;.\]

If we use these two expression in $- \ederiv \star \ederiv \star \kdifformh{\omega}{2} = \kdifform{f}{2}$, then the first operation we need to apply is the Hodge-$\star$ operator. This gives
\[ \star \kdifform{\omega}{2} = \star \sum_{i=1}^N \sum_{j=1}^N \omega_{ij}e_i(\xi)e_j(\eta) = \sum_{i=1}^N \sum_{j=1}^N \omega_{ij}\epsilon_i(\xi)\epsilon_j(\eta) = \sum_{i=1}^N \sum_{j=1}^N \tilde{\omega}_{ij}\tilde{h}_i(\xi) \tilde{h}_j(\eta) \;.\]
Here we used (\ref{eq:edge_functions_e_and_epsilon}) to write $e_i(\xi) = \epsilon(\xi)\ederiv \xi$, we use that $\star \ederiv \xi \ederiv \eta = 1$ and that we can represent this $0$-form exactly on the dual grid, since we have as many points on the dual grid as surfaces on the primal grid, by construction, see Figure~\ref{fig:mesh}. Whenever we refer to the dual grid we put twiddles, $\tilde{\cdot}$, on the coefficients and the basis forms.

\begin{remark}
Let the $2$-form be given by $\kdifform{\omega}{2}= \omega(\xi,\eta)\,\ederiv \xi \ederiv \eta$, then in 2D $\star \kdifform{\omega}{2}= \omega(\xi,\eta)$, because for $(\xi,\eta)\in [-1,1]^2$, the metric tensor is the identity. So the 'function' $\omega(\xi,\eta)$ remains exactly the same. But only $\ederiv \xi \ederiv \eta$ is chopped away from the expression. In that sense, the $\star$-operator can be seen as an identity operator -- because it does not change the functions $\omega(\xi,\eta)$ -- but only associates this 'function' with a different geometric object. This is generally not true in physical space $\manifold{M}_m$, where the function $\omega(x,y)$ will change when the Hodge is applied. The pullback therefore provides the metric connection. In the above case, initially the $2$-form was associated to surface and after the operation it is associated to points. What is less obvious is that the {\em pseudo-form} $\kdifform{\omega}{2}$ is converted into a {\em true form} $\tilde{\omega}^{(0)}$. The same happens for the finite dimensional form $\kdifformh{\omega}{2}$.
\end{remark}

After the $\star$-operation we are on the dual grid and there we want to apply the exterior derivative. This would give us
\begin{equation} 
\ederiv \sum_{i=0}^{N+1} \sum_{j=0}^{N+1} \tilde{\omega}_{ij}\tilde{h}_i(\xi) \tilde{h}_j(\eta) = \sum_{i=1}^{N+1} \sum_{j=1}^N
( \tilde{\omega}_{i,j} - \tilde{\omega}_{i-1,j}) \tilde{e}_i(\xi) \tilde{h}_j(\eta) + \sum_{i=1}^{N} \sum_{j=1}^{N+1}
( \tilde{\omega}_{i,j} - \tilde{\omega}_{i,j-1}) \tilde{h}_i(\xi) \tilde{e}_j(\eta) \;.
\label{eq:Hodge_of_zero_form}
\end{equation}
Note that here we refer to the values $\tilde{\omega}_{0,j}$ and $\tilde{\omega}_{N+1,j}$, for $j=1,\dots,N$ and to $\tilde{\omega}_{i,0}$ and $\tilde{\omega}_{i,N+1}$, for $i=1,\dots,N$. That is, we refer to values in points which are not part of the dual grid, see top right plot in Figure~\ref{fig:mesh}. It was already remarked, that the dual grid is not a complex and we need to add $0$-cells and $1$-cells to turn this grid into a cell complex. These additional cells are depicted on the bottom right in Figure~\ref{fig:mesh}. If these additional points constitute points on the boundary of the domain and if Dirichlet boundary conditions are prescribed, then these additional unknowns in (\ref{eq:Hodge_of_zero_form}) have a known value. If not, we simply add these boundary points to the scheme to have a cell complex for the dual grid. If the element depicted in Figure~\ref{fig:mesh} is one of the many elements in the spectral element mesh, we use these additional points to connect the solution between elements. 

\begin{remark}\label{rem:Dirichlet_dual}
When homogeneous Dirichlet boundary conditions are employed, the differences $( \tilde{\omega}_{i,j} - \tilde{\omega}_{i-1,j})$ and $( \tilde{\omega}_{i,j} - \tilde{\omega}_{i,j-1})$ essentially constitute the action of the coboundary on the dual grid as discussed in Section~\ref{Section::AlgebraicTopology} on page \pageref{page:coboundary_incidence}. The incidence matrix on the dual grid is connected to the incidence matrix on the primal grid by $\tilde{\mathsf{E}}_{(n-k,n-k-1)} = \incidenceboundary{k}{k-1}^T$. And the incidence matrix on the primal grid was determined by the connectivity and orientation of the primal grid. So here we explicitly see that the topological structure of the mesh in fact determines the derivatives on both the primal and the dual grid.
\end{remark}

\begin{remark}
In general, the approximation of $n$-forms, $\kdifform{\omega}{n}$, described in the way above, will lead to a discontinuous solution between elements. In fact, the trace of $n$-forms is not defined. The connectivity is obtained through $\star \kdifform{\omega}{n}$. In staggered finite volume methods, the introduction of these additional boundary points is also common and they are usually referred to as {\em 'ghost points'}. In the single grid approach, boundary values also need to be added which leads to the boundary integrals in weak formulations, see Remark~\ref{rem:Dirichlet_single}. The language differs, but the geometric structure is identical.
\end{remark}

Next we apply the $\star$-operator to (\ref{eq:Hodge_of_zero_form}) to obtain an expansion of the form
\begin{equation} 
\star \ederiv \star \kdifformh{\omega}{2} = - \sum_{i=1}^N \sum_{j=0}^N q_{i,j}^\eta e_i(\xi) h_j(\eta) + \sum_{i=0}^N \sum_{j=1}^N q_{i,j}^\xi h_i(\xi) e_j(\eta) \;.
\label{eq:fluxes_on_primal_grid}
\end{equation}
The Hodge expresses the two components in (\ref{eq:Hodge_of_zero_form}) on the dual grid in exactly the same components on the primal grid. So in that sense, the $\star$-operator is just a change of basis. The main thing, however, is that it converts the {\em true form} (\ref{eq:Hodge_of_zero_form}) into the {\em pseudo-form} (\ref{eq:fluxes_on_primal_grid}).

Finally, we take the exterior derivative of (\ref{eq:fluxes_on_primal_grid}) to obtain
\[ \ederiv \star \ederiv \star \kdifformh{\omega}{2} = \sum_{i=1}^N\sum_{j=1}^N \left ( q_{i,j}^\xi - q_{i-1,j}^\xi + q_{i,j}^\eta - q_{i,j-1}^\eta \right ) e_i(\xi) e_j(\eta) \;.\]
For the Poisson equation for a $2$-form, we need to equate this to $\kdifformh{f}{2}$ which gives
\[ \sum_{i=1}^N\sum_{j=1}^N \left ( q_{i,j}^\xi - q_{i-1,j}^\xi + q_{i,j}^\eta - q_{i,j-1}^\eta - f_{i,j} \right ) e_i(\xi) e_j(\eta) \;.\]
Since the basis functions $e_i(\xi)$ are all linearly independent, the only way in which we can satisfy this equation is by setting the expansion coefficients to zero, so
\[ q_{i,j}^\xi - q_{i-1,j}^\xi + q_{i,j}^\eta - q_{i,j-1}^\eta = f_{i,j} \;.\]
Given the explicit location of the dual grids, we can set up a matrix representation of the $\star$-operators in which case the Poisson equation for the $2$-form on dual grids can be written as
\begin{equation}
\mathsf{E}_{(2,1)} \tilde{H}^{1,1} \mathsf{E}_{(2,1)}^T H^{0,2} \psi(\kdifformh{\omega}{2}) = \psi(\kdifformh{f}{2}) \;.
\label{eq:discrete_system_dual}
\end{equation}
Here the incidence matrix $\mathsf{E}_{(2,1)}$ is the incidence matrix of the oriented primal cell complex, which depends on the connectivity of the grid and is therefore purely geometric. The incidence matrices are independent of the basis functions. The Hodge matrices do depend on the basis functions, which is another way of saying that the Hodge is a metric-dependent operator. The map $\psi$ is the map discussed on page \pageref{page:isomorphism} and given by (\ref{eq:isomorphism_cochains_expansion_coeff}). If we deform the grid, the incidence matrices will remain the same and only the metric-dependent part, i.e. the Hodge matrices will change as discussed in \citep{bouman::icosahom2009,jasper::eccomas2010}.

\subsection{Single grid approach}
The dual grid approach uses two grids and an explicit construction of the Hodge matrix. In the single grid approach the action of the Hodge operator is incorporated implicitly by using an inner product. The exterior derivative can be exactly represented on a grid due to the commutation relation (\ref{eq:commutation_relations_projection_derivative}). But we do not have such a relation for the codifferential operator. In the single grid approach, we therefore want to convert codifferentials into exterior derivatives using (\ref{eq:integration_by_parts}). 
In order to do so, we need to rewrite the equation $\ederiv \ederiv^* \kdifform{\omega}{2}=\kdifform{f}{2}$ in an equivalent first order system given by
\begin{equation}
\ederiv \ederiv^* \kdifform{\omega}{2}=\kdifform{f}{2} \quad \Longleftrightarrow \quad \left \{ \begin{array}{l}
\kdifform{q}{1} - \ederiv^* \kdifform{\omega}{2} = 0 \\
 \\
\ederiv \kdifform{q}{1} = \kdifform{f}{2}
\end{array} \right .
\end{equation}
If we take the inner product of the first equation with an arbitrary $1$-form, $\kdifform{v}{1}$ and the inner product of the second equation with an arbitrary $2$-form $\kdifform{p}{2}$, we obtain
\begin{equation} 
\left \{ \begin{array}{l}
\left ( \kdifform{q}{1},\kdifform{v}{1}\right )_{L^2\Lambda^1} - \left ( \ederiv^* \kdifform{\omega}{2},\kdifform{v}{1}\right )_{L^2\Lambda^1} = 0 \\
 \\
\left (\ederiv \kdifform{q}{1},\kdifform{p}{2}\right )_{L^2\Lambda^1} = \left ( \kdifform{f}{2},\kdifform{p}{2}\right )_{L^2\Lambda^1}
\end{array} \right .
\end{equation}
Now we can apply (\ref{eq:integration_by_parts}) to the second term in the first equation to obtain: Find $\kdifform{q}{1} \in H\Lambda^1$ and $\kdifform{\omega}{2} \in L^2\Lambda^2$ such that $\forall \kdifform{v}{1} \in H\Lambda^1$ and $\forall \kdifform{p}{2} \in L^2\Lambda^2$ we have
\begin{equation} 
\left \{ \begin{array}{l}
\left ( \kdifform{q}{1},\kdifform{v}{1}\right )_{L^2\Lambda^1} - \left ( \kdifform{\omega}{2},\ederiv \kdifform{v}{1}\right )_{L^2\Lambda^2} + \int_{\partial \manifold{M}}\mbox{tr}\,\kdifform{v}{1}\wedge \mbox{tr}\,\star \kdifform{\omega}{2}= 0 \\
 \\
\left (\ederiv \kdifform{q}{1},\kdifform{p}{2}\right )_{L^2\Lambda^1} = \left ( \kdifform{f}{2},\kdifform{p}{2}\right )_{L^2\Lambda^1}
\end{array} \right .
\label{eq:mixed_formulation}
\end{equation}
This formulation is equivalent to the mixed formulation in finite element methods, \citep{brezzi1991mixed}.

\begin{remark}
Where in the dual grid approach we needed to impose that the forms were {\em locally integrable}, we need to put more strict conditions on the admissible solution for the single grid approach. The solution must belong to suitably chosen Sobolev spaces.
\end{remark}

\begin{remark}\label{rem:Dirichlet_single}
Dirichlet boundary conditions imposed on $\kdifform{\omega}{2}$ appear through the boundary integral (weak imposition of boundary conditions). Here we see that we need to prescribe $\mbox{tr}\,\star \kdifform{\omega}{2}$, just as in the case for the dual grid approach, see Remark~\ref{rem:Dirichlet_dual}. If we prescribe $\mbox{tr}\,\star \kdifform{\omega}{2}$, the boundary integral can be transferred to the righthand side. If Neumann conditions are prescribed, i.e. $\mbox{tr}\,\kdifform{q}{1}$ is given, then $\mbox{tr}\,\kdifform{v}{1}=0$ and the contribution of the boundary integral vanishes.
\end{remark}

Well-posedness of the weak formulation can be found in any book on mixed finite element methods, for instance, \citep{brezzi1991mixed}. Since the exterior derivative maps $\Lambda^1$ onto $\Lambda^2$, we can always find a particular $\kdifform{q}{1}_f$ such that $\ederiv \kdifform{q}{1}_f = \kdifform{f}{2}$. Let ${\mathcal Z} = \{ \kdifform{v}{1} \in H\Lambda^1 \,|\, \ederiv \kdifform{v}{1} = 0 \,\}$, then if we restrict the solution space and the space of weight forms to ${\mathcal Z}$ we obtain: Find $\kdifform{q}{1} \in {\mathcal Z}$, such that for all $\kdifform{v}{1} \in {\mathcal Z}$, we have 
\[ \left ( \kdifform{q}{1},\kdifform{v}{1} \right )_{L^2\Lambda^1} = -\int_{\partial \manifold{M}}\mbox{tr}\,\kdifform{v}{1}\wedge \mbox{tr}\,\star \kdifform{\omega}{2} - \left ( \kdifform{q}{1}_f,\kdifform{v}{1} \right )_{L^2\Lambda^1} \;.\]
The bilinear form $\left ( \kdifform{q}{1},\kdifform{v}{1} \right )_{L^2\Lambda}$ is trivially coercive and bounded. Consider the space ${\mathcal Z}^\perp=\{ \kdifform{a}{1}\in H\Lambda^1\,|\, (\kdifform{a}{1},\kdifform{v}{1})_{L^2\manifold{M}}\,, \forall \kdifform{v}{1} \in {\mathcal Z}\,\}$. Since the map $\ederiv\,:\,{\mathcal Z}^\perp \rightarrow L^2 \Lambda^2$ is a bijection for $n=2$, we have the Poincar\'{e} inequality that for every $\kdifform{p}{2} \in L^2\Lambda^2$ there exists a unique $\kdifform{v}{1}_p \in {\mathcal Z}^\perp \subset H\Lambda^1$ such that $\ederiv \kdifform{v}{1}_p = \kdifform{p}{2}$ and there exists a constant $c_p>0$ such that, \citep{arnold2006finite,arnold2010finite,Kreeft_A_priori_Stokes}
\[ \| \kdifform{v}{1}_p \|_{L^2\Lambda^2} \leq c_p \| \kdifform{p}{2} \|_{L^2\Lambda^1}\;.\]
Then it follows that, \citep{Kreeft_A_priori_Stokes,Bochev_VKI}
\[ \sup_{\kdifform{v}{1} \in L^2\Lambda^2} \frac{\left ( \kdifform{p}{2},\ederiv \kdifform{v}{1}\right )_{L^2\Lambda^2}}{\|\kdifform{v}{1}\|_{L^2\Lambda^2}} \geq \frac{\left ( \kdifform{p}{2},\ederiv \kdifform{v}{1}_p\right )_{L^2\Lambda^2}}{\|\kdifform{v}{1}_q\|_{L^2\Lambda^2}} = \frac{\| \kdifform{p}{2}\|^2_{L^2\Lambda^2}}{\|\kdifform{v}{1}_p\|_{L^2\Lambda^2}} \geq \frac{1}{c_p} \| \kdifform{p}{2}\|^2_{L^2\Lambda^2} \;.\]
Therefore, the inf-sup condition, \citep{brezzi1991mixed}, is satisfied at the continuous level.

If we restrict the infinite dimensional function spaces to finite dimensional conforming subspaces using the mimetic projection $\projection$, we obtain the discrete equation: Find $\kdifformh{q}{1} \in H\Lambda_h^1$ and $\kdifformh{\omega}{2} \in L^2\Lambda_h^2$ such that $\forall \kdifformh{v}{1} \in H\Lambda_h^1$ and $\forall \kdifformh{p}{2} \in L^2\Lambda_h^2$ we have
\begin{equation} 
\left \{ \begin{array}{l}
\left ( \kdifformh{q}{1},\kdifformh{v}{1}\right )_{L^2\Lambda^1} - \left ( \kdifformh{\omega}{2},\ederiv \kdifformh{v}{1}\right )_{L^2\Lambda^2} + \int_{\partial \manifold{M}}\mbox{tr}\,\kdifformh{v}{1}\wedge \mbox{tr}\,\star \kdifformh{\omega}{2}= 0 \\
 \\
\left (\ederiv \kdifformh{q}{1},\kdifformh{p}{2}\right )_{L^2\Lambda^1} = \left ( \kdifformh{f}{2},\kdifformh{p}{2}\right )_{L^2\Lambda^1}
\end{array} \right .
\label{eq:discrete_mixed_formulation}
\end{equation}
Well-posedness of this discrete formulation follows directly from the fact that $H\Lambda_h^1 \subset H\Lambda^1$, $L^2\Lambda_h^2 \subset L^2\Lambda^2$ and ${\mathcal Z}_h = \{ \kdifformh{v}{1} \in H\Lambda_h^1 \,|\, \ederiv \kdifformh{v}{1} = 0 \, \} \subset {\mathcal Z}$. The latter property is a direct consequence of (\ref{eq:commutation_relations_projection_derivative}).

\vskip 0.3cm
\noindent
\framebox[1.0\linewidth]{
\begin{minipage}{0.92\linewidth}
If in the continuous setting, the equation $\ederiv \kdifform{q}{1} = \kdifform{f}{2}$ possesses a solution $\kdifform{q}{1}$, then the discrete solution $\ederiv \kdifformh{q}{1} = \kdifformh{f}{2}$ has a solution $\projection \kdifform{q}{1}$, \citep{kreeft2011mimetic}
\[ 0 = \projection \left ( \kdifform{f}{2} - \ederiv \kdifform{q}{1} \right ) \stackrel{(\ref{eq:commutation_relations_projection_derivative})}{=} \projection \kdifform{f}{2} - \ederiv \projection \kdifform{q}{1} = \kdifformh{f}{2} - \ederiv \kdifformh{q}{1} \;.\]
\end{minipage}}
\vskip 0.3cm 

\begin{remark}
The inf-sup condition is not really a 'condition', because it is identically satisfied by the mimetic projection. For Stokes flow this is proven in \citep{Kreeft_A_priori_Stokes}.
\end{remark}

The discrete matrix equation corresponding to (\ref{eq:discrete_mixed_formulation}) has the following form
\begin{equation}
\left(
\begin{array}{cc}
-\matrixoperator{M}^1 & \left ( \matrixoperator{M}^{2}\incidenceboundary{2}{1} \right )^T \\
\matrixoperator{M}^2 \incidenceboundary{2}{1} & \matrixoperator{0}
\end{array}
\right )
\left (
\begin{array}{c}
\psi(\kdifformh{q}{1}) \\
\psi(\kdifformh{\omega}{2})
\end{array}
\right)
=
\left(
\begin{array}{c}
\matrixoperator{\mbox{tr}\,\star \kdifformh{\omega}{2}} \\
\psi(\kdifformh{f}{2})
\end{array}
\right)\;,
\label{eq:single_grid_matrix_system}
\end{equation}
in which $\matrixoperator{M}^1$ is the mass matrix associated with the $1$-form basis forms, (\ref{eq:basis_k_forms}) and $\matrixoperator{M}^2$ is the mass matrix for the $2$-forms, (\ref{eq:basis_k_forms}). The incidence matrix, $\incidenceboundary{2}{1}$, denotes differentiation at cochain level and is again determined by the geometry. The map $\psi$ is the isomorphism which connects the cochain to its expansion coefficients, (\ref{eq:isomorphism_cochains_expansion_coeff}). Elimination of $\psi\kdifformh{q}{1})$ from this system yields
\[ \matrixoperator{M}^2 \incidenceboundary{2}{1} \left ( \matrixoperator{M}^1 \right )^{-1} \incidenceboundary{2}{1} \left ( \matrixoperator{M}^2 \right ) \psi(\kdifformh{\omega}{2}) = \matrixoperator{M}^2 \psi(\kdifformh{f}{2}) \;,\]
or after pre-multiplication by the inverse of $\matrixoperator{M}^2$
\begin{equation}
\incidenceboundary{2}{1} \left ( \matrixoperator{M}^1 \right )^{-1} \incidenceboundary{2}{1} \left ( \matrixoperator{M}^2 \right ) \psi(\kdifformh{\omega}{2}) = \psi(\kdifformh{f}{2}) \;.
\label{eq:discrete_system_single}
\end{equation}
If we compare this with the discrete equation for the dual grid method, (\ref{eq:discrete_system_dual}), we note that in the single grid method the mass matrices and their inverses play the role of the Hodge matrix. This relation between mass matrices and Hodge matrices was also observed by \citep{bossavit1999::japan_02,bochev2006principles,arnold2006finite,hiptmair2001discrete,tarhasaari1999some}.


\subsection{Results for the Poisson problem for  volume forms}

\begin{figure}[ht]
\centering
\subfigure{
\includegraphics[width=0.25\textwidth]{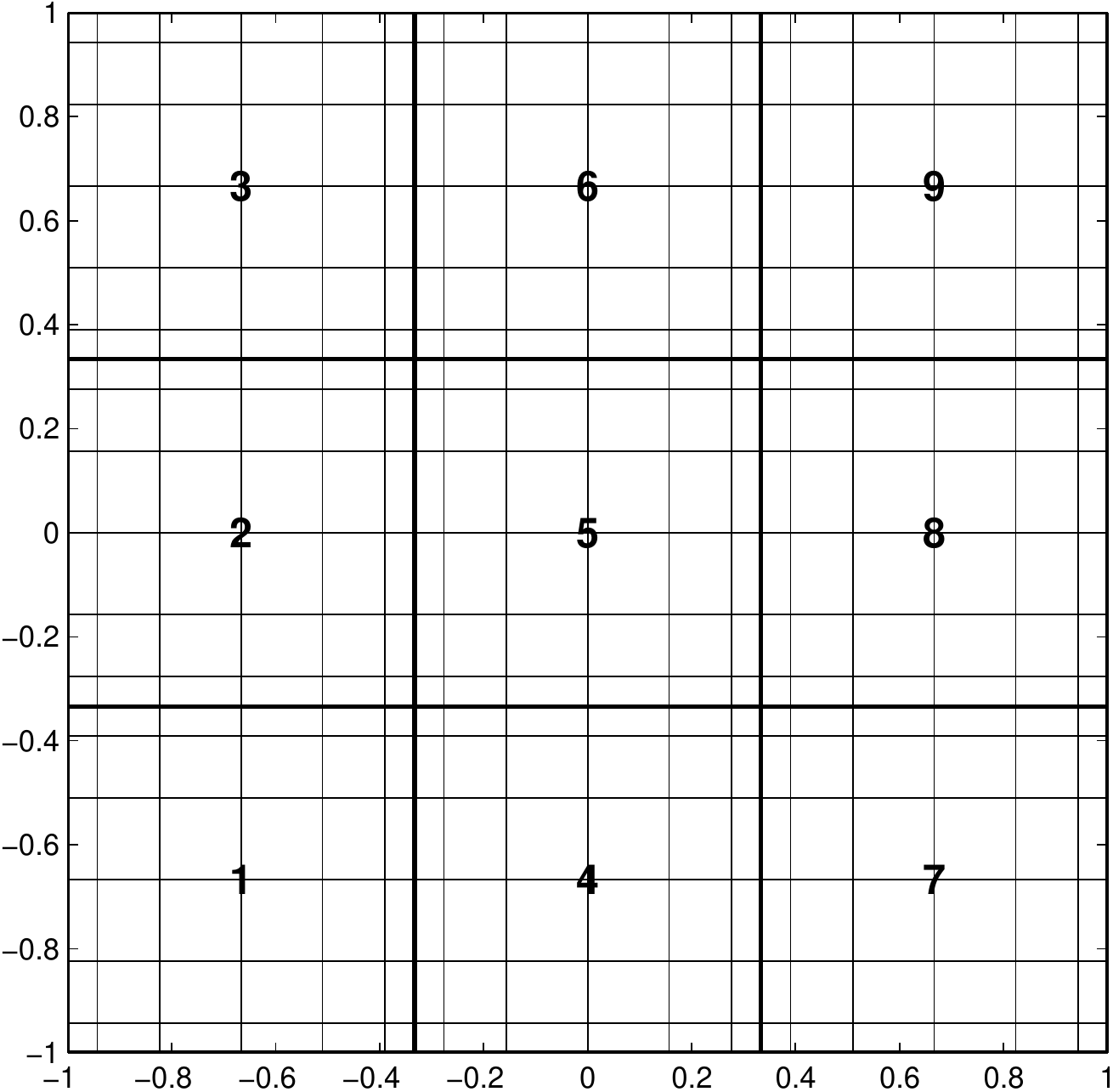}
\label{fig:crazyMesh_cc000}
}
\subfigure{
\includegraphics[width=0.25\textwidth]{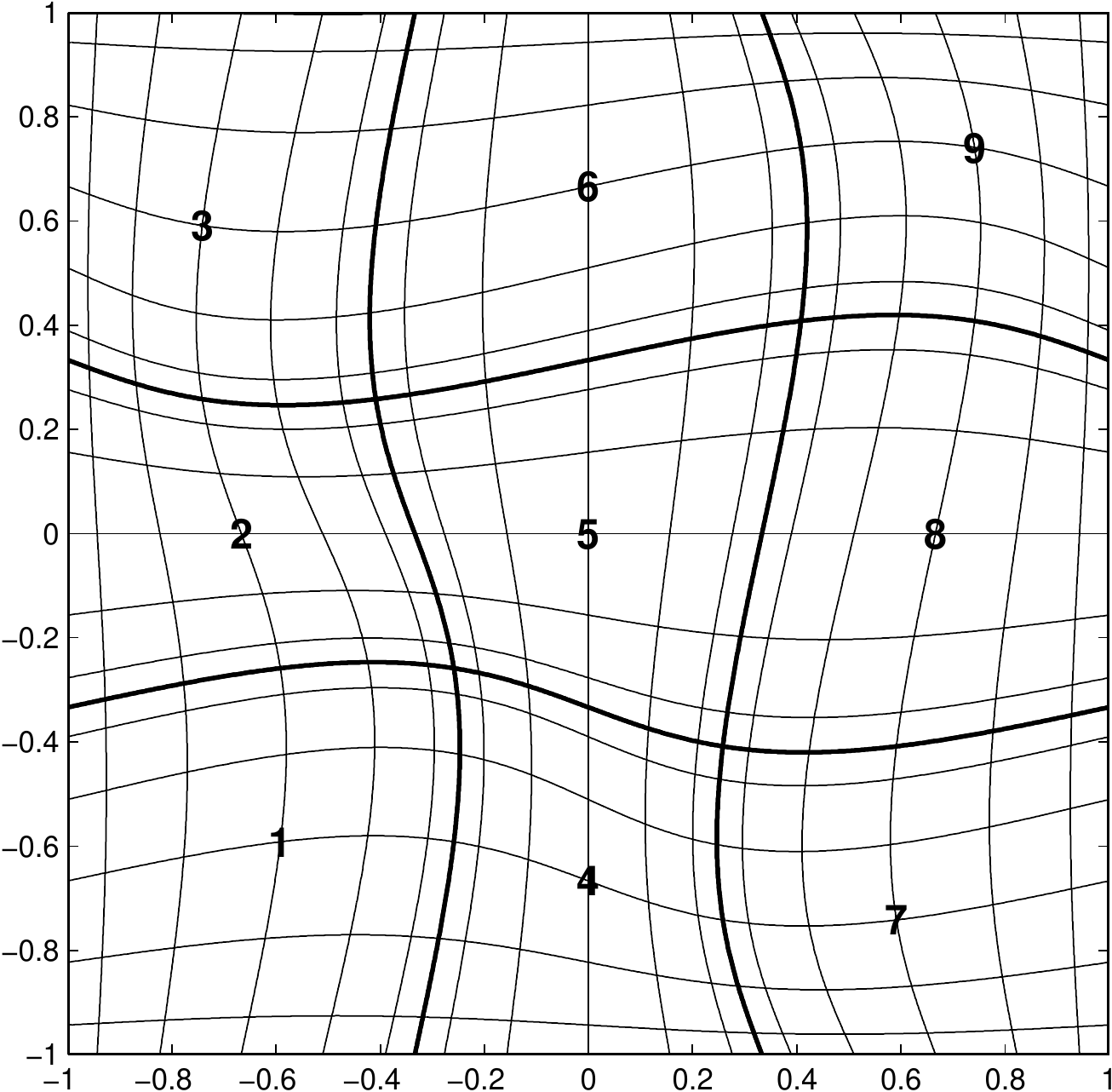}
\label{fig:crazyMesh_cc010}
}
\subfigure{
\includegraphics[width=0.25\textwidth]{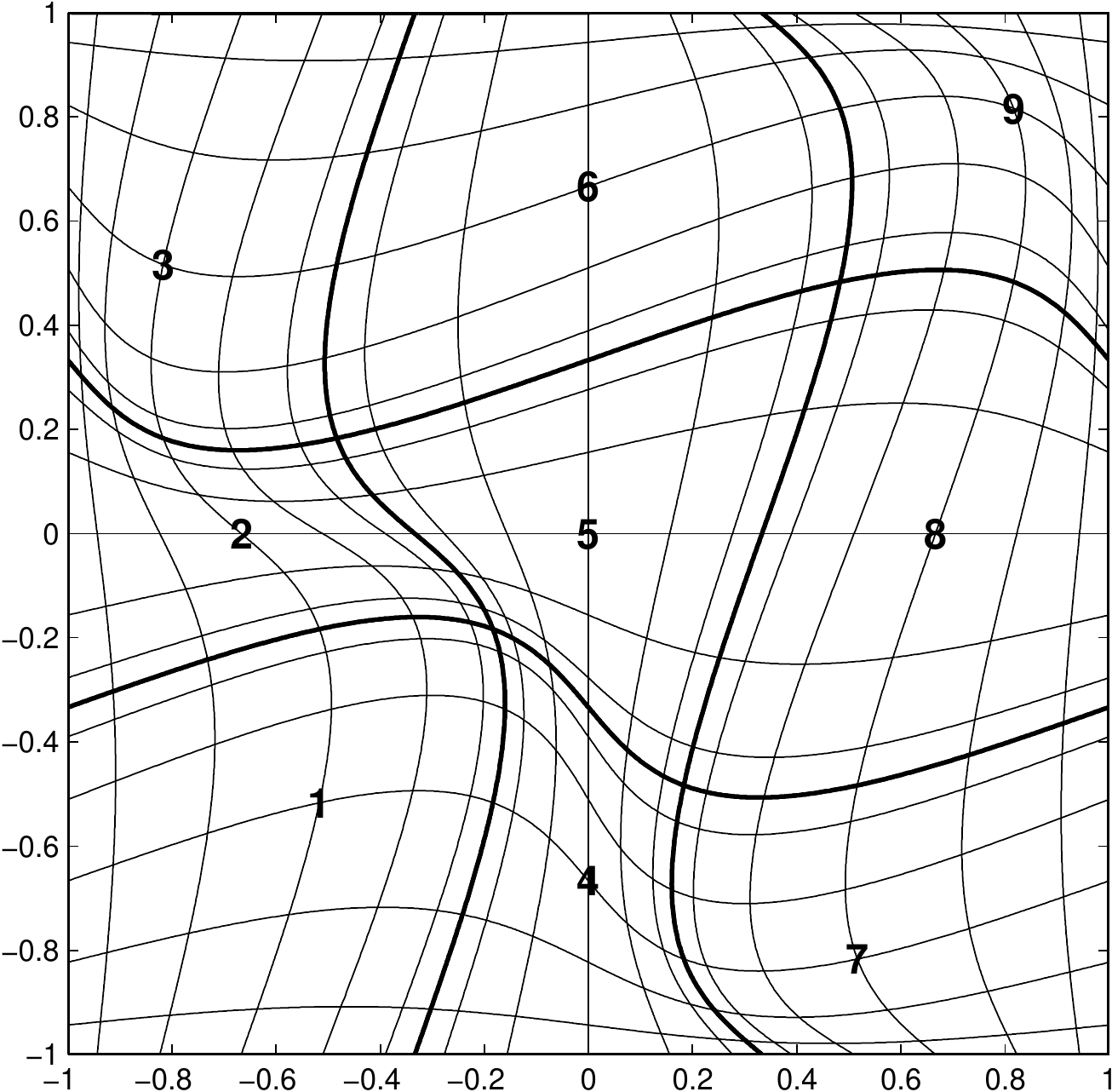}
\label{fig:crazyMesh_cc020}
}
\caption{Deformed geometry mesh for deformation coefficient $c=0$ (left), $c=0.1$ (middle) and $c=0.2$ (right), for $3\times3$ elements of order $N=6$.}
\label{fig:crazyMesh_cc010_cc020}
\end{figure}

In this section we will present results for a sample problem using both the dual grid and the single grid. As already remarked in Example~\ref{ex:incidence_matrix}, the incidence matrices remain invariant under distortion of the grid. And (\ref{eq:discrete_system_dual}) and (\ref{eq:discrete_system_single}) reveal that differentiation is determined by these incidence matrices. When we deform the grid, only the Hodge matrices in (\ref{eq:discrete_system_dual}) and the mass matrices in (\ref{eq:discrete_system_single}) will change, \citep{bouman::icosahom2009,palha::icosahom2009}.
Therefore, for this test case we consider three meshes obtained by a transformation of the unit domain $[-1,1]^2$ to curvilinear coordinates. The mapping is given by $(x,y)=\Phi(\xi,\eta)$, with
\begin{subequations}
\label{mapping1}
\begin{align}
x(\xi,\eta)&=\xi+c\sin(\pi\xi)\sin(\pi\eta)\\
y(\xi,\eta)&=\eta+c\sin(\pi\xi)\sin(\pi\eta).
\end{align}
\end{subequations}
For $c=0$ we obtain the orthogonal Gauss-Lobatto grid. The grid obtained for $c=0.0$, $c=0.1$ and $c=0.2$ are displayed in Figure~\ref{fig:crazyMesh_cc010_cc020}

the exact solution for this problem is given by
\begin{equation}
\kdifform{\alpha}{2}=\sin(2\pi x)\sin(2\pi y)\,\ederiv x \ederiv y \;.
\end{equation}
For these tests Dirichlet boundary conditions corresponding to the exact solution are prescribed.

\subsubsection{Results for dual grid approach}
\begin{figure}[ht!]
\centering
\subfigure{
\includegraphics[width=0.3\textwidth]{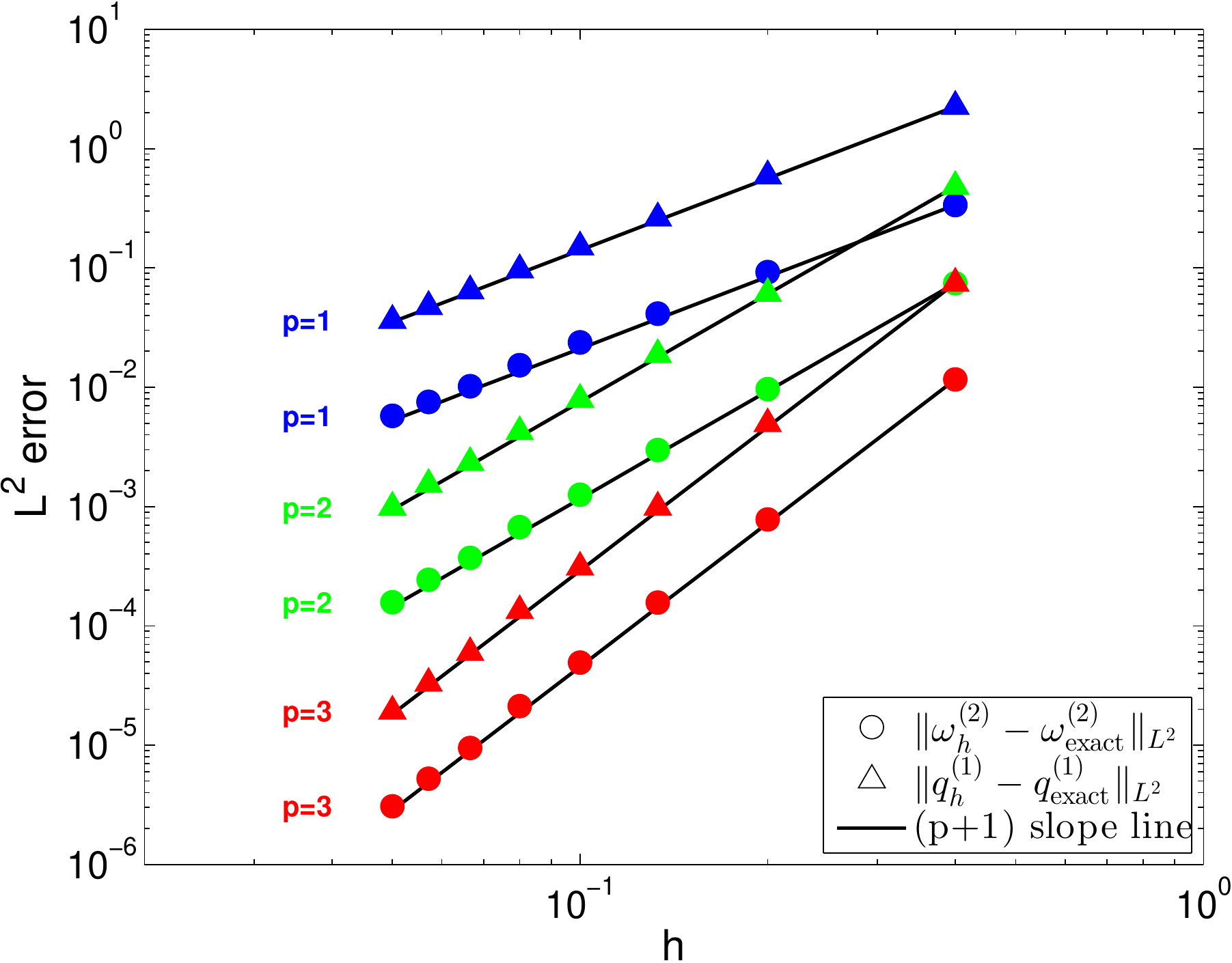}
\label{fig:h_crazy_000_convergence_primal_dual}
}
\subfigure{
\includegraphics[width=0.3\textwidth]{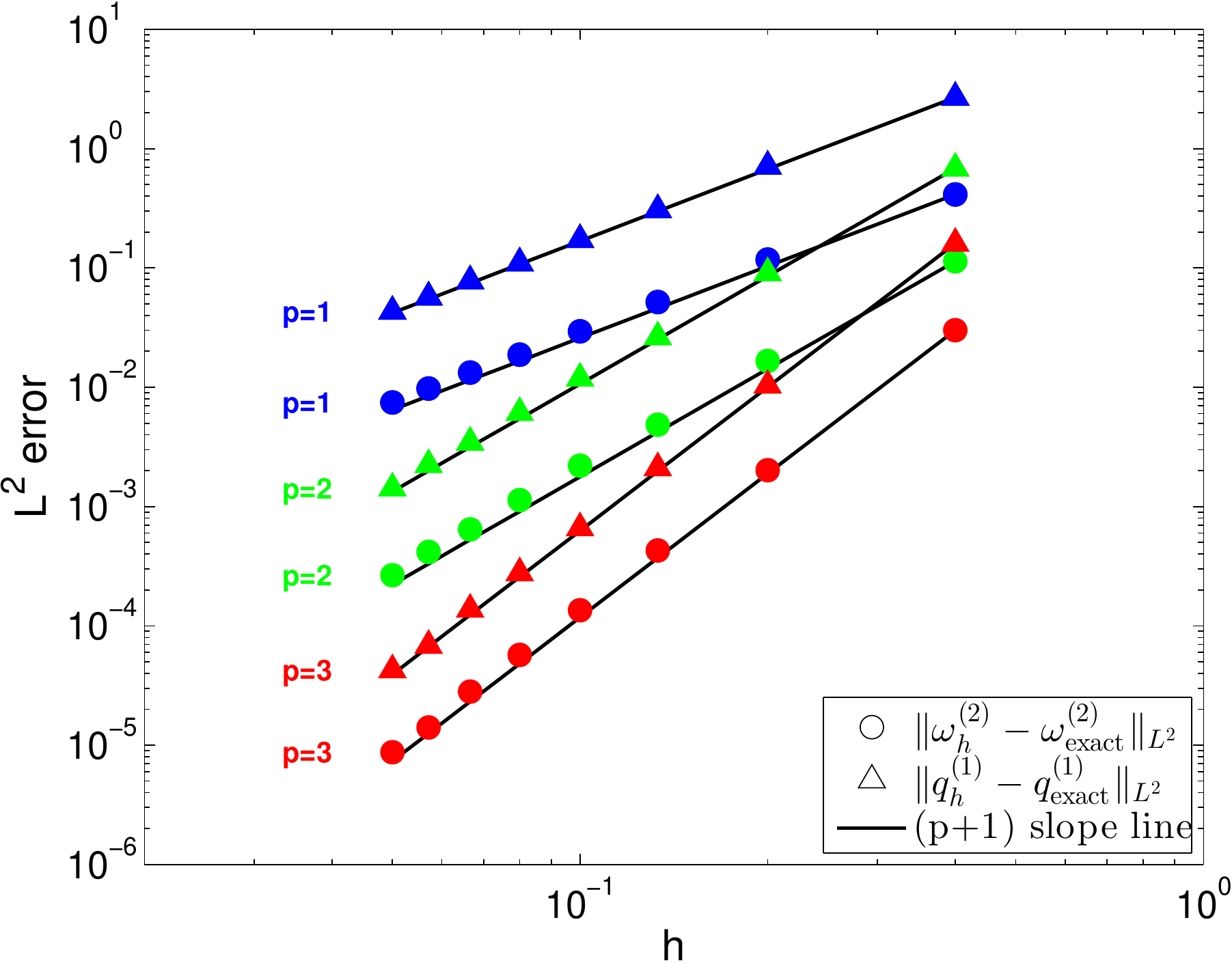}
\label{fig:h_crazy_010_convergence_primal_dual}
}
\subfigure{
\includegraphics[width=0.3\textwidth]{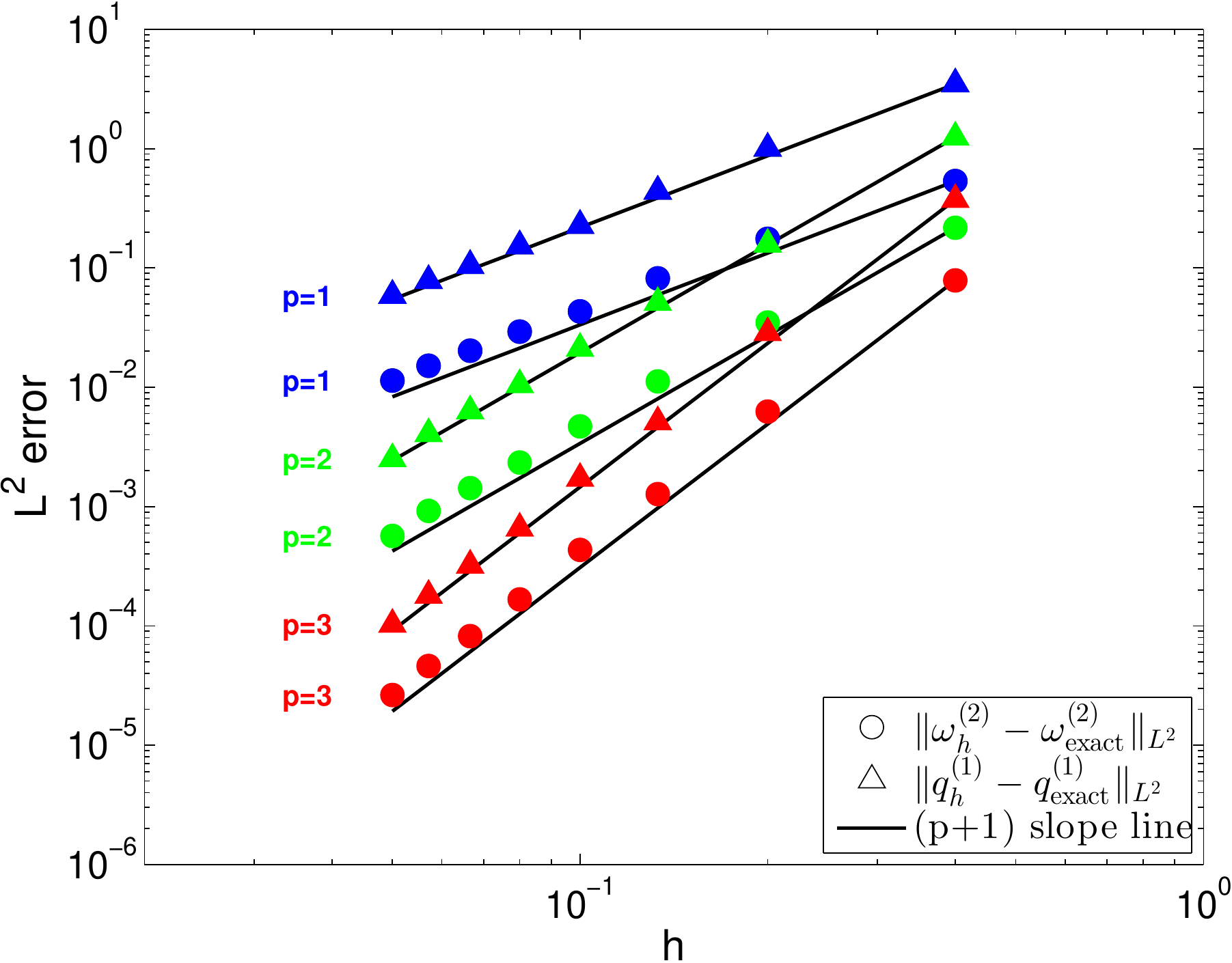}
\label{fig:h_crazy_020_convergence_primal_dual}
}
\caption{Plots of $h$-convergence of the $L^{2}$ error for $N=1,2,3$ on the meshes shown in Figure~\ref{fig:crazyMesh_cc010_cc020} for the dual grid approach \eqref{eq:discrete_system_dual} and for the single grid approach \eqref{eq:single_grid_matrix_system}. Left plot mesh $c=0.0$, middle plot $c=0.1$ and right plot for $c=0.2$. The straight lines indicate the slope of the expected convergence rate.}
\label{fig:h_crazy_convergence_primal_dual}
\end{figure}

In Figure~\ref{fig:h_crazy_convergence_primal_dual} $h$-convergence is shown for polynomial degree $N=1,2,3$, respectively, for the meshes shown in Figure~\ref{fig:crazyMesh_cc010_cc020}.

The $h$-convergence curves for the dual grid approach and the single grid approach are plotted in Figure~\ref{fig:h_crazy_convergence_primal_dual} and are indistinguishable. The expected rate of convergence for both the solution $\kdifformh{\omega}{2}$ and the fluxes $\kdifformh{q}{1}$ are expected to be $N+1$, where $N$ is the polynomial degree, which is confirmed by the results shown in Figure~\ref{fig:h_crazy_convergence_primal_dual}.
\begin{figure}[ht!]
\centering
\subfigure{
\includegraphics[width=0.3\textwidth]{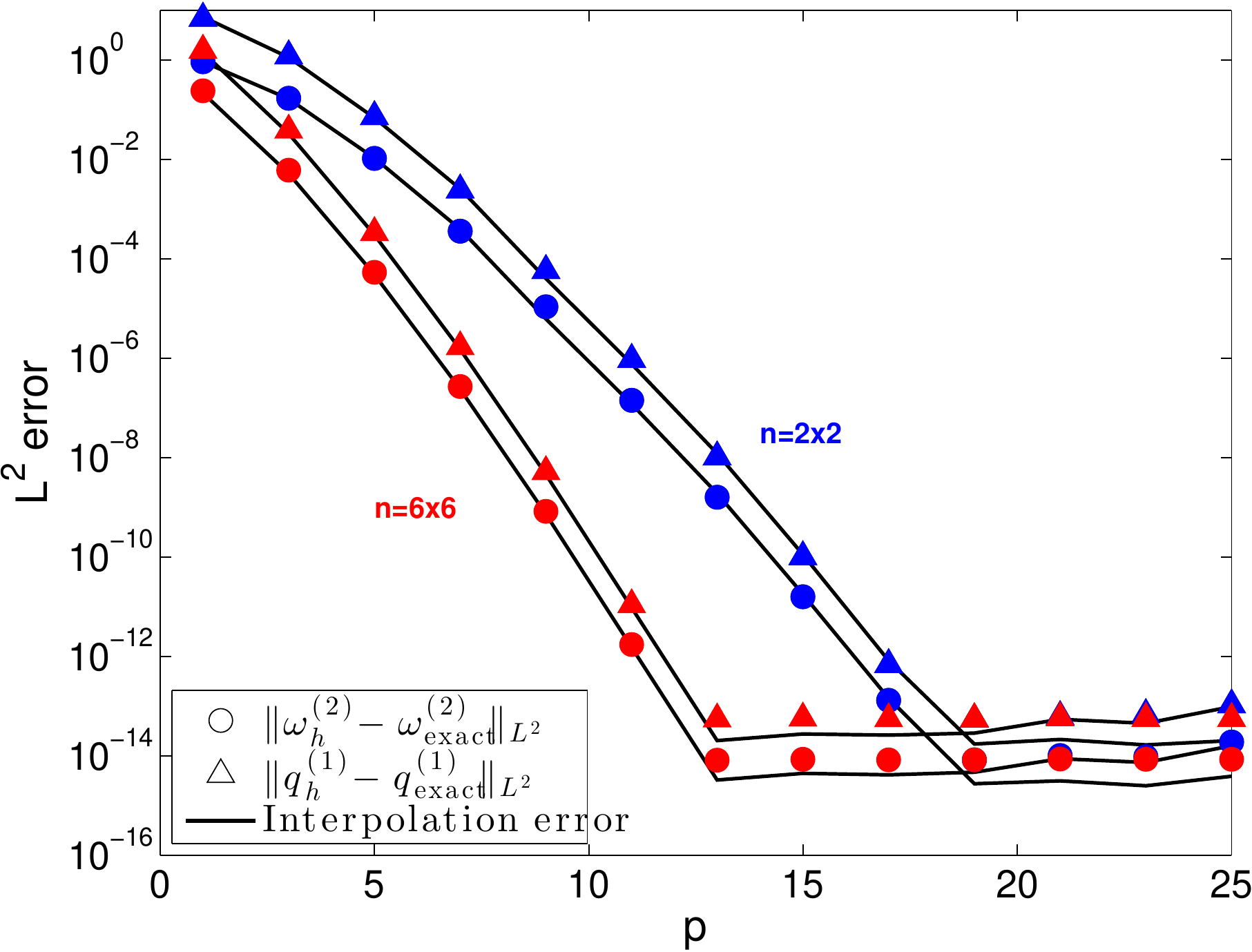}
\label{fig:p_crazy_000_convergence_primal_dual}
}
\subfigure{
\includegraphics[width=0.3\textwidth]{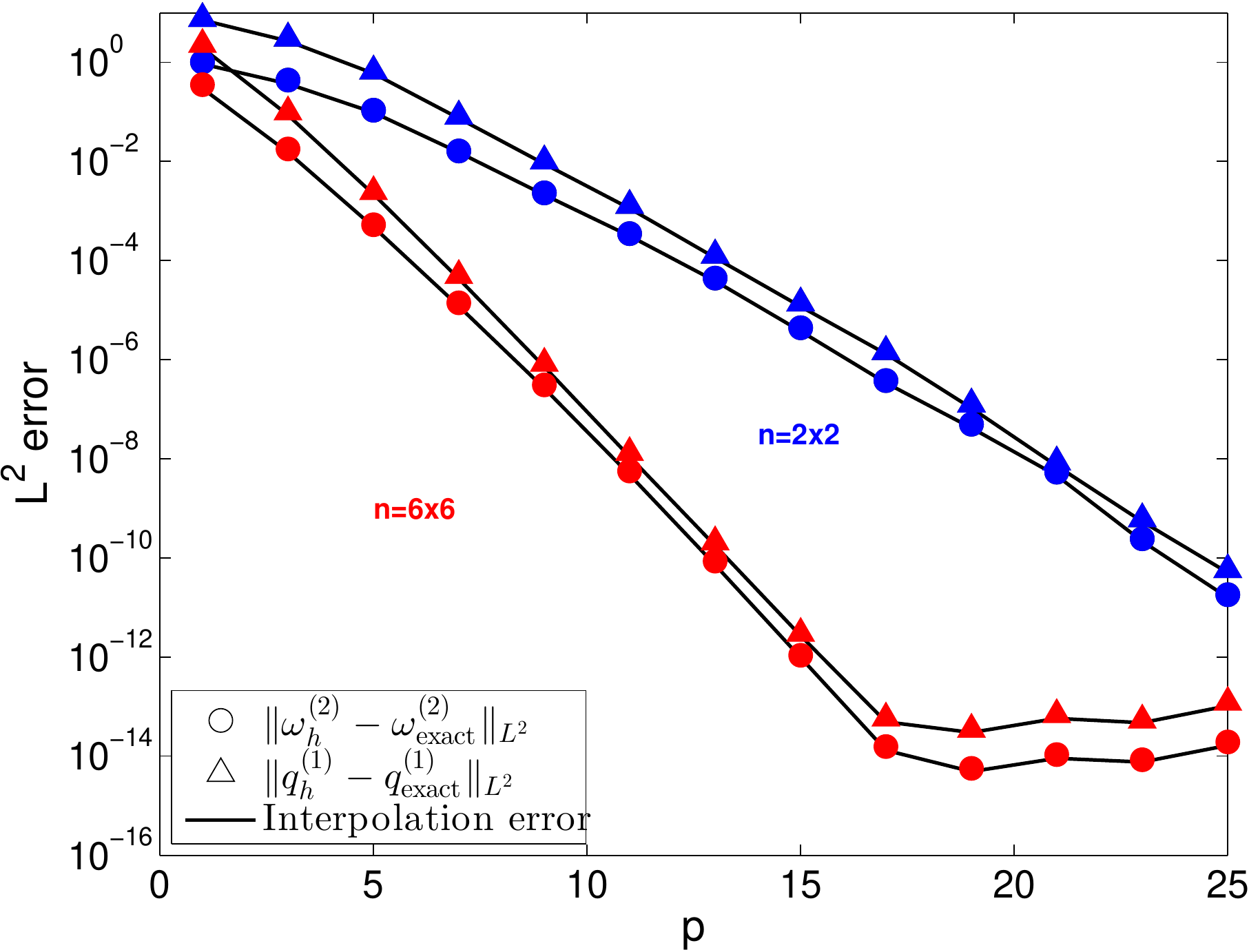}
\label{fig:p_crazy_010_convergence_primal_dual}
}
\subfigure{
\includegraphics[width=0.3\textwidth]{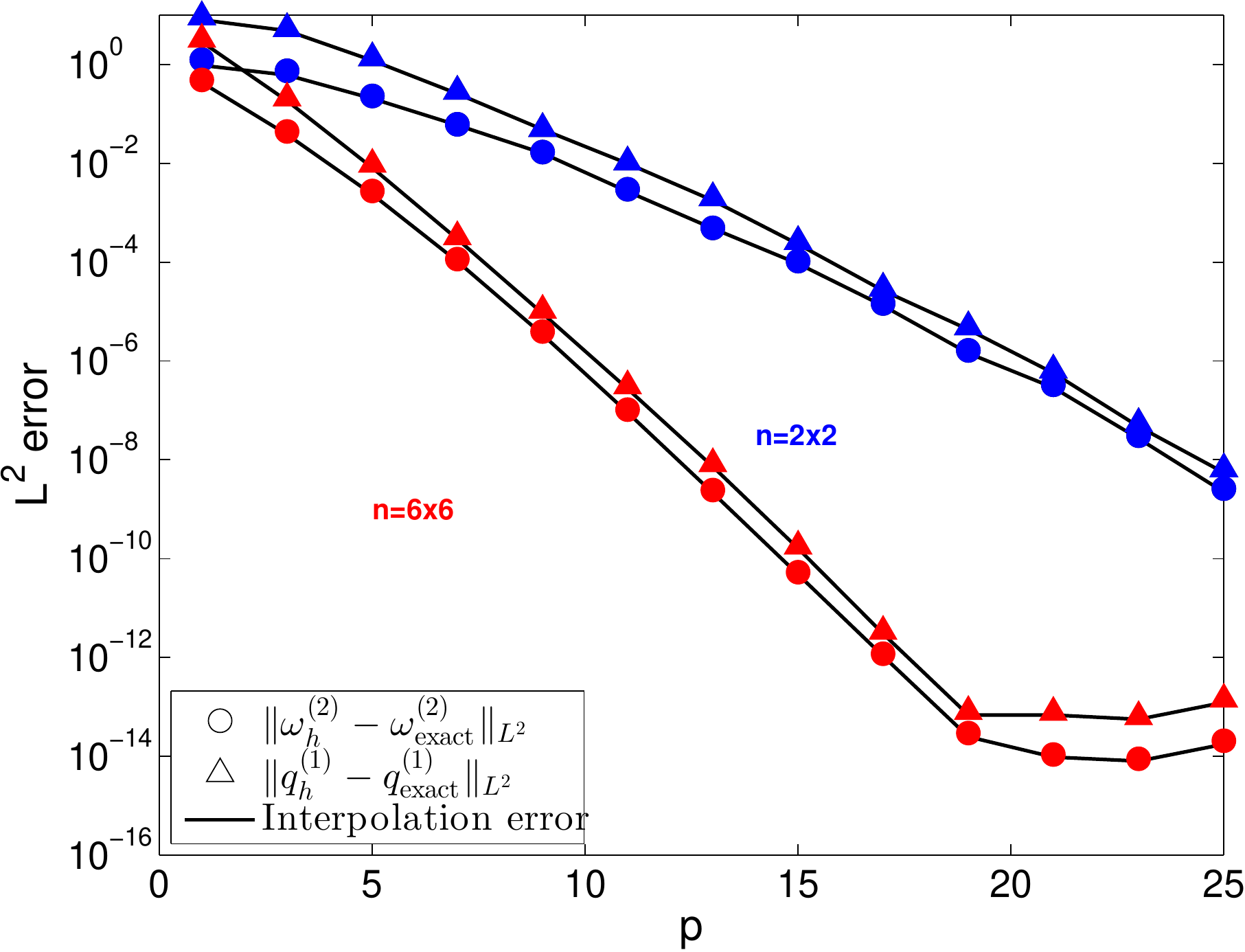}
\label{fig:p_crazy_020_convergence_primal_dual}
}
\caption{Plots of $p$-convergence of the $L^{2}$ error on a $2\times 2$ mesh and a $4\times 4$ mesh on the meshes shown in Figure~\ref{fig:crazyMesh_cc010_cc020} for the dual grid approach \eqref{eq:discrete_system_dual}and for the single grid approach \eqref{eq:single_grid_matrix_system}. Left plot mesh $c=0.0$, middle plot $c=0.1$ and right plot for $c=0.2$. The black line is the interpolation of the exact solution.}
\label{fig:p_crazy_convergence_primal_dual}
\end{figure}
In Figure~\ref{fig:p_crazy_convergence_primal_dual} $p$-convergence is plotted for both the dual grid approach and the single grid method. Both methods display identical convergence rates. The solid line in Figure~\ref{fig:p_crazy_convergence_primal_dual} is the $L^2$-error of the interpolation of the exact solution.
Despite the fact that the convergence rates are exactly the same for both methods, the solutions are not identical.

\begin{figure}[ht!]
\centering
\subfigure{
\includegraphics[width=0.35\textwidth]{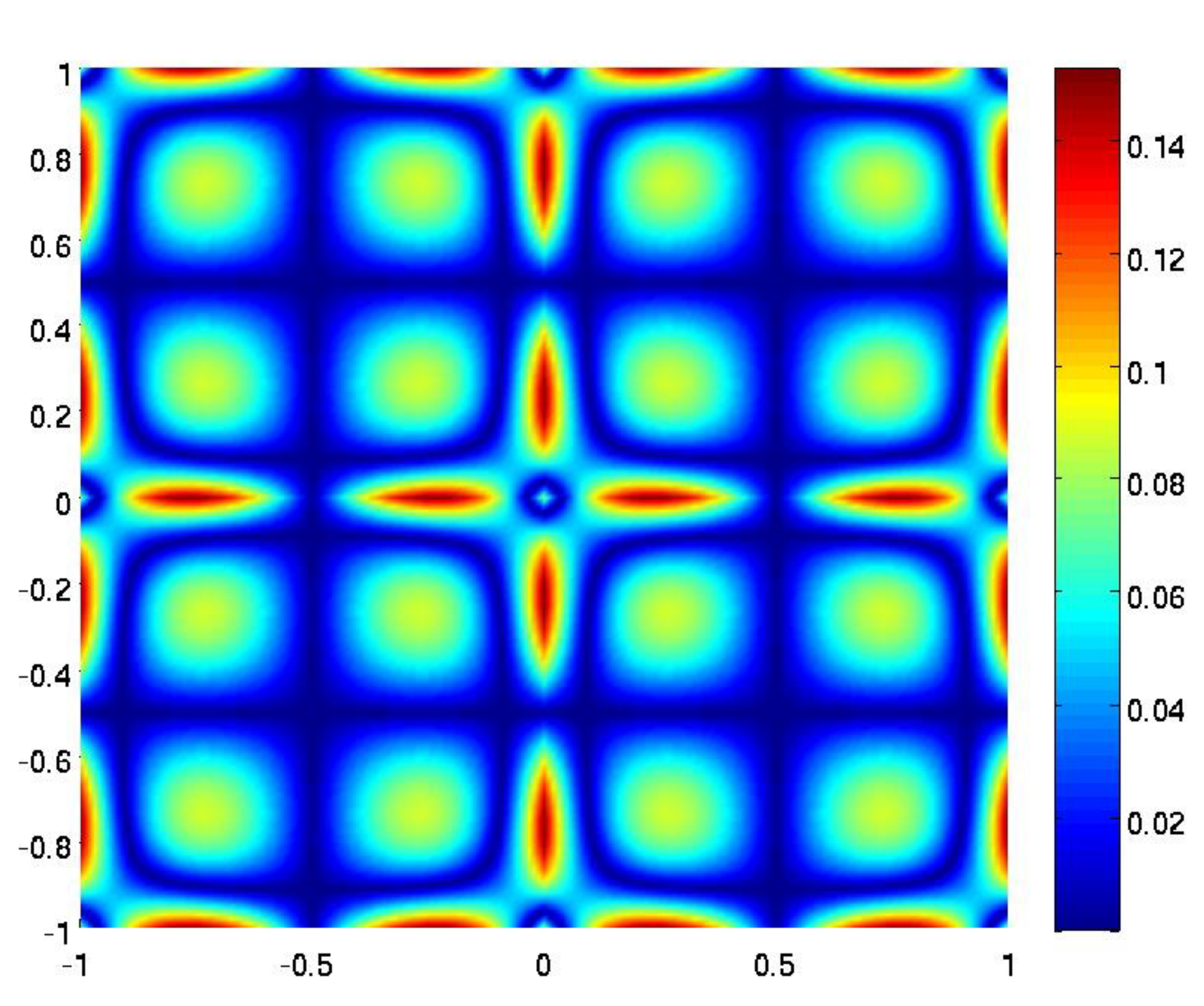}
\label{fig:diff_dual_single_p3}
}
\subfigure{
\includegraphics[width=0.35\textwidth]{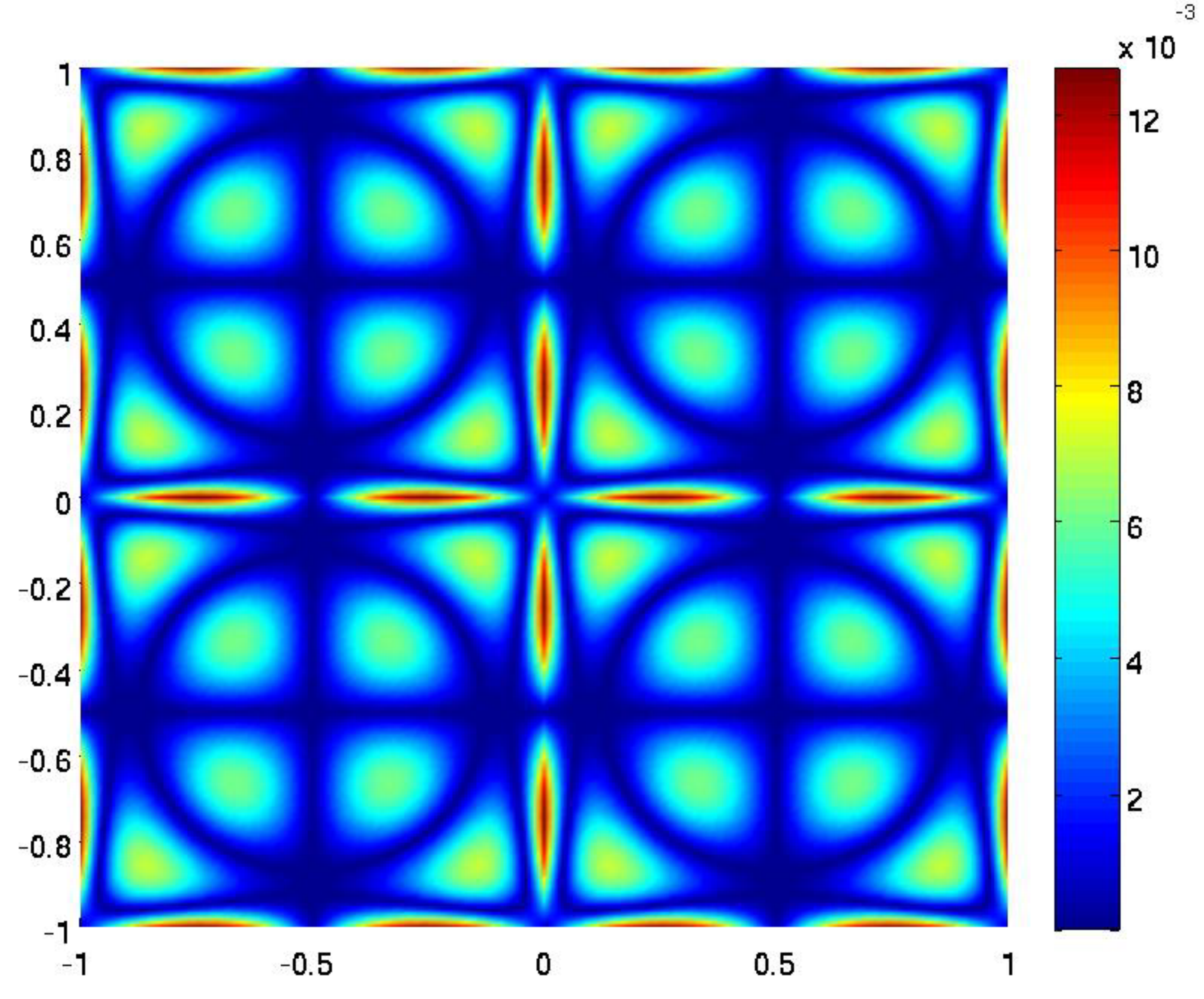}
\label{fig:diff_dual_single_p5}
}
\caption{Contour plot of the absolute value of the spatial difference between the dual grid solution and the single grid solution on a $2 \times 2$ grid, with $p=3$ (left) and $p=5$ (right).}
\label{fig:diff_dual_single}
\end{figure}
\begin{figure}[ht!]
\centering
\subfigure{
\includegraphics[width=0.3\textwidth]{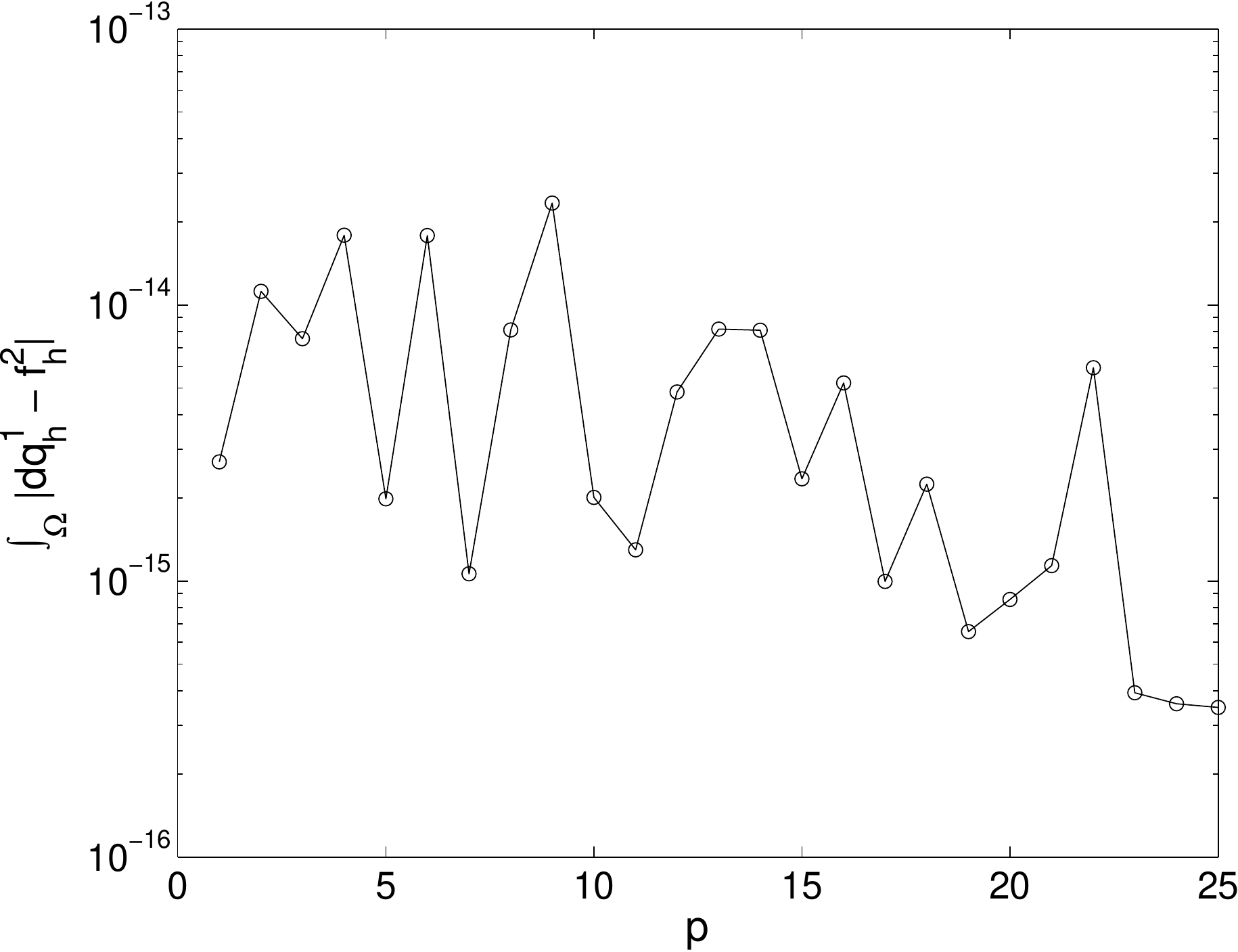}
\label{fig:div_q=f_c00_dual}
}
\subfigure{
\includegraphics[width=0.3\textwidth]{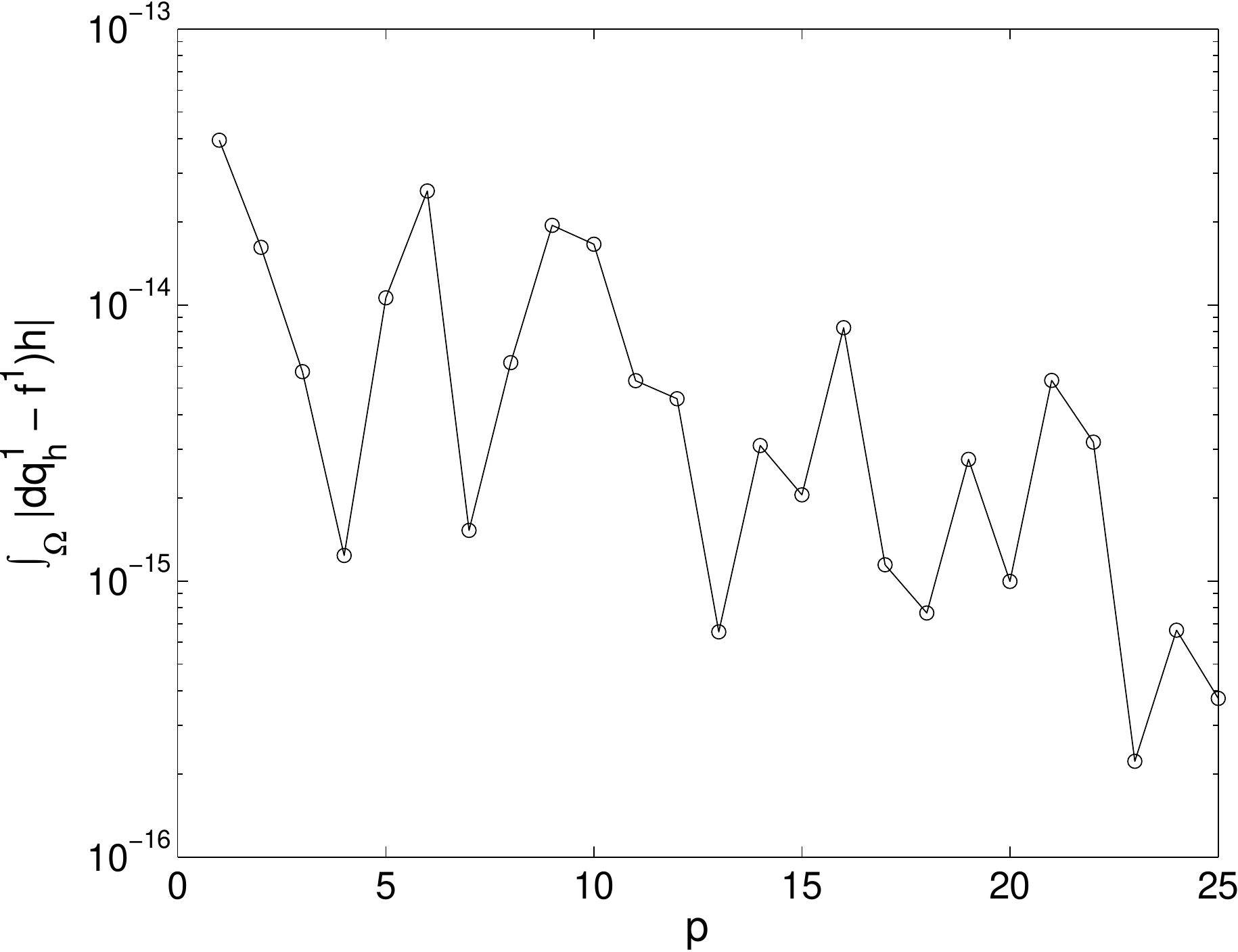}
\label{fig:div_q=f_c01_dual}
}
\subfigure{
\includegraphics[width=0.3\textwidth]{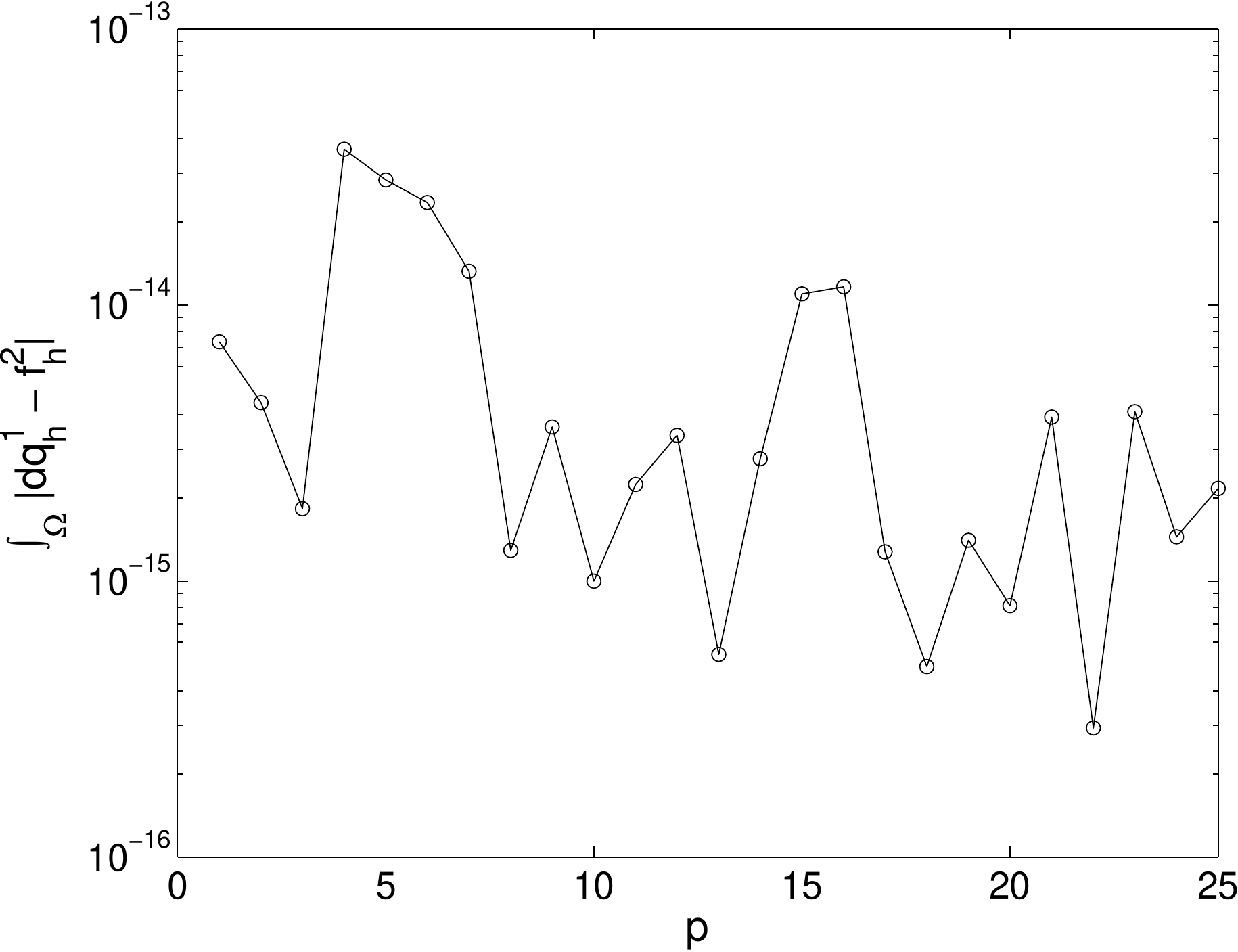}
\label{fig:div_q=f_c02_dual}
}
\caption{$L^\infty$-error for the \underline{dual grid approach} in the equation $\ederiv \kdifformh{q}{1}=\kdifformh{f}{2}$ for $c=0.0$ (left), $c=0.1$ (middle) and $c=0.2$ (right) }
\label{fig:div_q=f_dual}
\end{figure}
\begin{figure}[ht!]
\centering
\subfigure{
\includegraphics[width=0.3\textwidth]{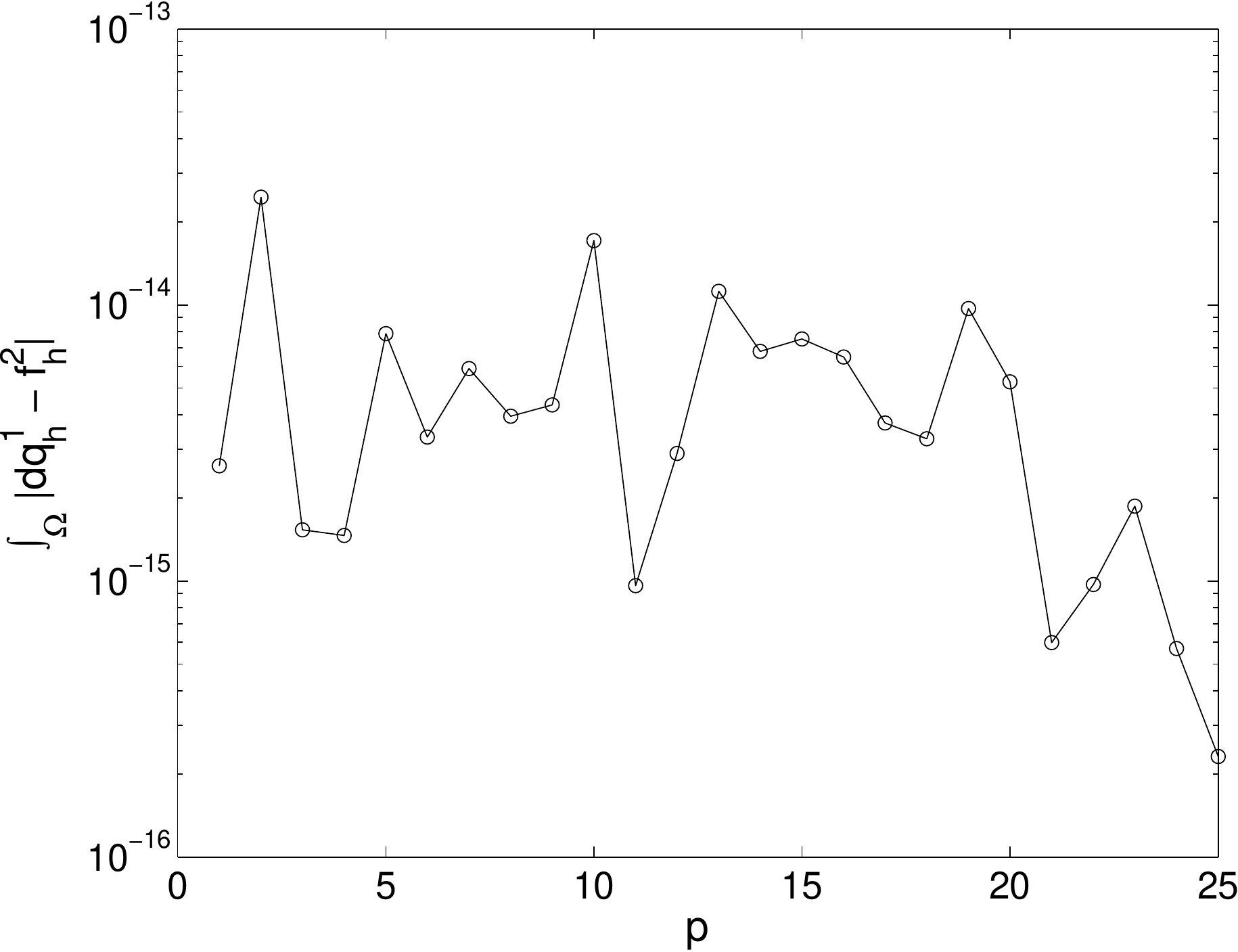}
\label{fig:div_q=f_c00_single}
}
\subfigure{
\includegraphics[width=0.3\textwidth]{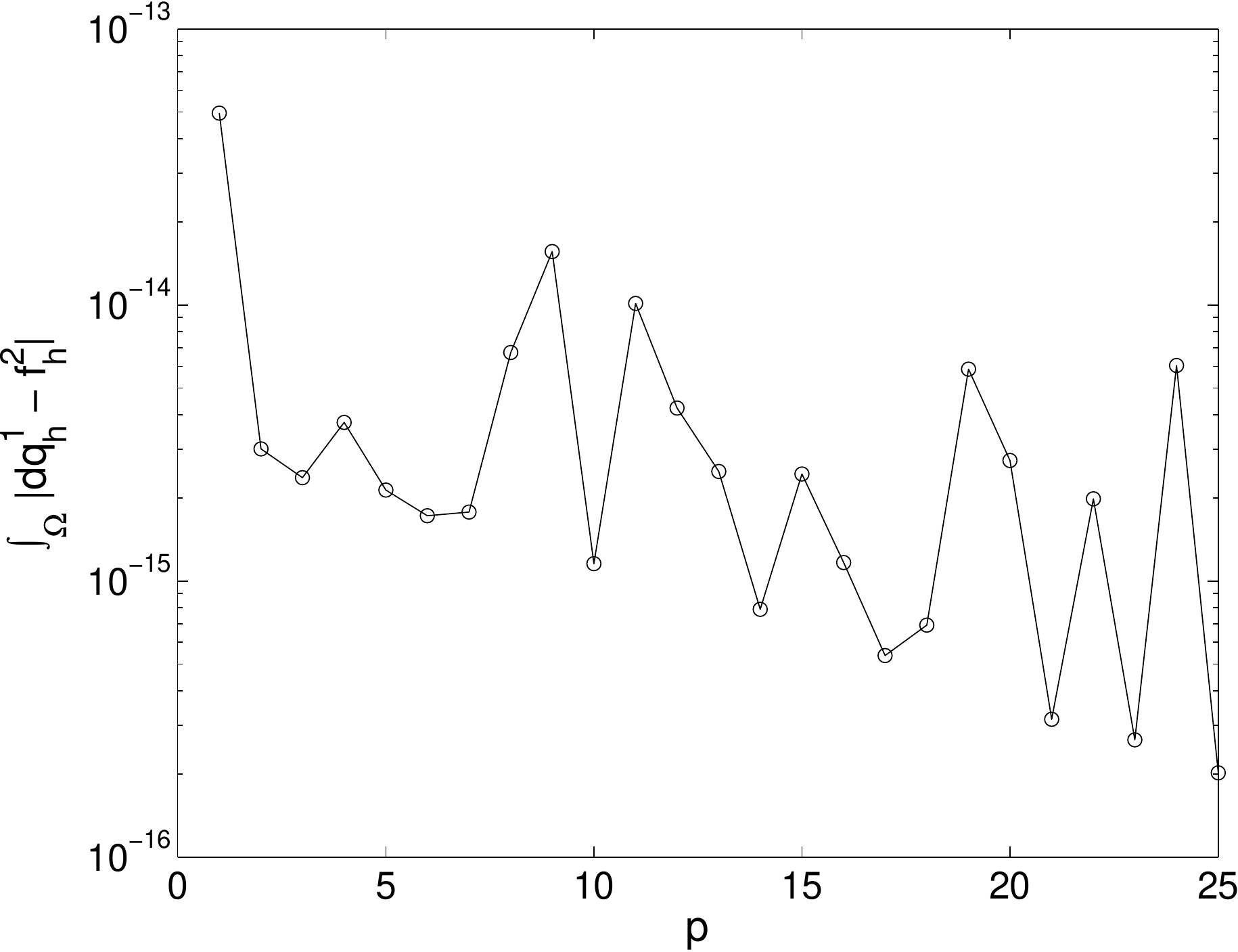}
\label{fig:div_q=f_c01_single}
}
\subfigure{
\includegraphics[width=0.3\textwidth]{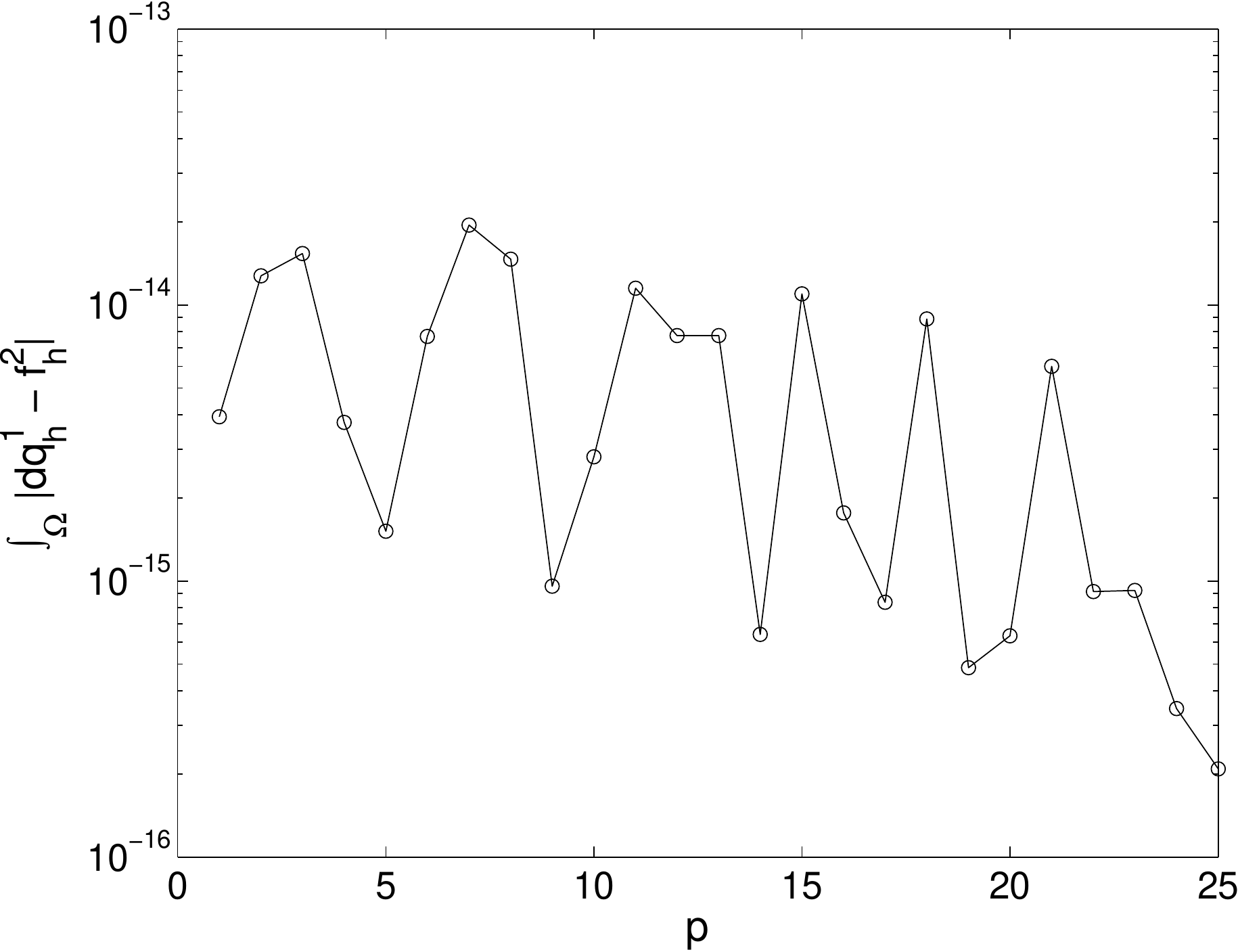}
\label{fig:div_q=f_c02_single}
}
\caption{$L^\infty$-error for the \underline{single grid approach} in the equation $\ederiv \kdifformh{q}{1}=\kdifformh{f}{2}$ for $c=0.0$ (left), $c=0.1$ (middle) and $c=0.2$ (right) }
\label{fig:div_q=f_single}
\end{figure}

Figure~\ref{fig:diff_dual_single} displays the absolute value of the spatial difference between the solution obtained with a dual grid approach and the single grid approach. While the fluxes $\kdifformh{q}{1}$ in both methods will be slightly different, the conservation law $\ederiv \kdifformh{q}{1} = \kdifformh{f}{2}$ is satisfied for both methods as shown in Figure~\ref{fig:div_q=f_dual} and Figure~\ref{fig:div_q=f_single}. For this particular test problem $\kdifform{f}{2} = -8\pi^2 \sin(2\pi x) \sin (2\pi y)\,\ederiv x \ederiv y$.

\section{SUMMARY AND OUTLOOK}
\label{Section::Conclusions}

In this paper we presented two formulations for the numerical solution of the Laplace equation: One in which the action of the Hodge operator is explicitly performed using two topological dual grids and one where the action of the Hodge is embedded in the definition of the inner product. In both methods the discrete derivative, which is the formal adjoint of the geometric boundary operator, is explicit in terms of the incidence matrix of oriented grid.

In terms of the $L^2$-convergence rate for the Poisson problem for volume forms, for $h$- and $p$-refinement, both methods are optimal and the convergence plots show indistinguishable curves on orthogonal and curved meshes. The actual solution, however, differ.
The major difference between the \emph{primal-dual grid formulation} and the \emph{single grid mixed formulation} lies in the discretization of the codifferential operator. Despite these difference, conservation equation $\ederiv \kdifformh{q}{1}=\kdifformh{f}{2}$ is satisfied up to machine precision.

Current work focuses on eigenvalue problems for the Laplace operator applied to $k$-forms and the discretization of the Lie derivative in order to model convective behaviour.

\bibliographystyle{plainnat}
\bibliography{./library}

\end{document}